\documentclass[10pt,journal,onecolumn]{IEEEtran}
\usepackage{bm, amssymb, amsmath, graphicx, mathtools, soul, xcolor, float, algorithm, algorithmic, siunitx, booktabs, subfigure}

\usepackage[english]{babel}
\addto\captionsenglish{}
\addto\captionsenglish{}

\newcommand{\B}[1]{{\bm #1}}
\newcommand{\T}{^{\mbox{\tiny T}}}
\newcommand{\dd}{\; \text{d}}
\newcommand{\tfc}{\emph{Theory of Functional Connections}}
\newcommand{\ce}{\emph{constrained expression}}
\newcommand{\ces}{\emph{constrained expressions}}

\begin{document}

\title{Least-squares solutions of boundary-value problems in hybrid systems}

\author{Hunter Johnston\footnote{PhD Student, Aerospace Engineering, Texas A\&M University, College Station, TX $77843$-$3141$, USA. E-mail: hunterjohnston$@$tamu.edu.} and Daniele Mortari\footnote{D. Mortari, Professor, Aerospace Engineering, Texas A\&M University, College Station, TX $77843$-$3141$, USA. E-mail: mortari$@$tamu.edu.}
}

\maketitle

\begin{abstract}

    This paper looks to apply the mathematical framework of the \tfc\ to the solution of boundary-value problems arising from hybrid systems. The \tfc\ is a technique to derive \emph{constrained expressions} which are analytical expressions with embedded constraints. These expressions are particularly suitable to transform a large class of constrained optimization problems into unconstrained problems. The initial and most useful application of this technique is in the solution of differential equations where the problem can be posed as an unconstrained optimization problem and solved with simple numerical techniques (i.e. least-squares).

    A hybrid system is simply a sequence of different differential equations. The approach developed in this work derives an analytical \ce\ for the entire range of a hybrid system, enforcing both the boundary conditions as well as the continuity conditions across the sequence of differential equations. This reduces the searched solution space of the hybrid system to only admissible solutions. The transformation allows for a least-squares solution of the sequence for linear differential equations and a iterative least-squares solution for nonlinear differential equations.
    
    Lastly, the general formulation for ``$n$'' segments is developed and validation is provided through numerical tests for three differential equation sequences: linear/linear, linear/nonlinear, and nonlinear/nonlinear. The accuracy level obtained are all at machine-error, which is consistent with the accuracy experienced in past studies on the application of the \tfc\ to solve single ordinary differential equations.
\end{abstract}

\section{Introduction}

Hybrid systems are dynamical systems governed by a time-sequence of differential equations (DE), which can be either linear or nonlinear. A simple example is a bouncing ball where the motion is described by a sudden DE variation at impact. Hybrid systems also arise in many control problems where the system is subject to a bang-bang discrete controls. These are a special case of hybrid systems called variable structure systems (VSS) and the study of the control of these systems is called variable structure control (VSC) \cite{VSS}.

Techniques to solve hybrid systems have existed since the 1960s utilizing multiple shooting methods \cite{shooting_TPBVP,shooting_bang,shooting_BVP,LineShoot}. In these approaches, the interval is divided over multiple sub-intervals and the boundary-value problem (BVP) is converted to multiple initial-value problems (IVP). The unknown boundary conditions are then solved by minimizing the DE residuals and the residuals of function and derivative continuities connecting all sub-intervals. In practice, root solving techniques (bisection, Newton's method, etc.) are used to minimize all residuals. In general, even when two subsequent linear DEs are connected, solutions based on shooting methods require an initial guess of the unknown parameters that is used to iterate until the solution is obtained. Note that the convergence is dictated by the initial guess \cite{shooting_converg} and it is not guaranteed. Regardless, studies have been conducted to quantify the error of these methods once an approximation is obtained \cite{shooting_err1,shooting_err2}.

Other techniques for solving these problems include the finite different method where the DE is approximated by a difference equation which converts the problem into a system of equations that are solved using linear algebra techniques. On the other hand, in finite element methods (collocation, Galerkin, etc.) \cite{FEM}, the problem is split into smaller parts called finite elements and simple approximated equations are used to model these elements. These elements are then assembled into a larger system of equations that model the entire problem. The major drawback the finite difference and finite element method is the number of subdivisions needed to capture large variations in the solution. 

Recently, an alternative approach to solving DEs was developed based on \emph{Theory of Functional Connections}\footnote{This theory, initially called ``Theory of Connections,'' has been renamed for two reasons. First, the name ``Theory of Connections'' already identifies a specific theory in differential geometry, and second, what this theory is actually doing is ``Functional Interpolation'' as it provides \emph{all} functions satisfying a set of constraints in term of function and any derivative in rectangular domains in $n$-dimensional spaces.} (TFC) \cite{U-ToC, LDE, NDE, M-ToC}. In general, TFC is a mathematical framework that transforms constrained problems into unconstrained problems by deriving \ces. \emph{Constrained expressions} are mathematical expressions that, using a free function $g (x)$, describe the family of all possible functions that satisfy a given problem's constraints. In this way, TFC reduces the search space to only those functions that satisfy the problem's constraints, thus transforming the original constrained problem into an unconstrained problem. The general expression to derive univariate \ces\ follows the form,
\begin{equation}\label{eq:gen_CE}
    y (x) = g (x) + \sum_{k=1}^n \eta_k \, p_k (x)
\end{equation}
where $p_k (x)$ are $n$ assigned linearly independent functions and $g (x)$ is a free function. The $\eta_k$ are coefficient functions that are derived by imposing the set of $n$ constraints. The constraints considered in the TFC method are any linear combination of the functions and/or derivatives evaluated at specified values of the variable $x$. The constraints of the DE are used to directly solve for the unknown coefficient functions, $\eta_k$. Once the vector of the unknown $\B{\eta}$ is determined, the constraints for the DE are satisfied for any possible $g (x)$. 

TFC was originally introduced as a means to satisfy point constraints, derivative constraints, and relative constraints \cite{U-ToC} . Eventually, TFC was extended to encompass infinite and integral constraints \cite{constraints}, and recently, extended to include multivariate boundary value constraints and multivariate, arbitrary-order derivative constraints \cite{M-ToC}. Additionally, Ref. \cite{TFC_selected} highlights the broader application of this technique outside of DEs. 

In the area of DEs, TFC has proven to be a powerful tool for solving linear \cite{LDE} and nonlinear \cite{NDE} ordinary DE (initial-, boundary-, and multi-value problems). By expressing $g(x)$ as a set of known basis functions (e.g., Fourier series or orthogonal polynomials, such as Legendre or Chebyshev polynomials) with unknown coefficients, linear DEs can be solved with simple numerical techniques like least-squares, while nonlinear problems necessitate an iterative least-squares approach to converge to the desired solution \cite{NDE}. The major advantage of TFC over other numerical integrators is: 1) speed (it is orders of magnitudes faster when compared with Runge-Kutta methods), 2) machine level accuracy in solutions, and 3) it provides a unified framework to solve all classes of ordinary differential equations (ODE). This means the same methodology is used to solve initial-, boundary-, and multi-value problems, with the only difference being the specific \ce.

While most ODEs can be solved with the theory developed in Refs. \cite{LDE,NDE}, BVPs arising from hybrid systems require additional framework. For these systems, a \ce\ for each segment of continuous dynamics must be derived and their continuity enforced. This paper details a method (based on the \tfc\ approach) that produces the \ces\ that automatically satisfy both DE constraints and continuity constraints. Thus, these problems can also be converted into an unconstrained optimization problem with embedded and guaranteed continuities. First, a summary of the TFC method to solve a simple BVP is provided that will be used as a building block for the theory to follow. Following this, two approaches to embed continuity constraints are developed and compared. Lastly, this approach and methodology is validated for several hybrid systems.

\section{The TFC Approach to Solve BVPs}

A simple example of the TFC is shown here for the convenience of the reader, but more details can be found in Refs. \cite{U-ToC,NDE}. Let us consider solving a second-order boundary value problem such that,
\begin{equation}\label{eq:gen_DE}
    F(x,y,y',y'') = 0 \quad \text{subject to: }\begin{cases} y(x_0) = y_0 \\ y(x_f) = y_f\end{cases}
\end{equation}
The \ce\ can be searched using Eq. (\ref{eq:gen_CE}) where $p_1 (x) = 1$ and $p_2 (x) = x$ (see Ref. \cite{U-ToC,NDE} for details), which leads to,
\begin{equation}\label{eq:ex_ce}
    y(x) = g(x) + \eta_1 + \eta_2 \, x.
\end{equation}
By applying the constraints, we obtain the following system of equations,
\begin{equation*}
    \begin{Bmatrix} y_0 - g_0 \\ y_f - g_f\end{Bmatrix} = \begin{bmatrix} 1 & x_0 \\ 1 & x_f\end{bmatrix} \begin{Bmatrix} \eta_1 \\ \eta_2 \end{Bmatrix},
\end{equation*}
and by inverting the matrix we can solve for the unknown $\eta$ values which are,
\begin{equation*}
    \begin{aligned}
    \eta_1 &= \dfrac{1}{x_f-x_0} \Big[x_f\Big(y_0 - g_0)\Big) - x_0\Big(y_f - g_f)\Big)\Big]\\
    \eta_2 &= \dfrac{1}{x_f-x_0} \Big[\Big(y_f - g_f\Big) - \Big(y_0 - g_0)\Big)\Big].
    \end{aligned}
\end{equation*}
These can be substituted into Eq. (\ref{eq:ex_ce}) to obtain the \ce,
\begin{equation}\label{eq:ex_ce_solved}
    y(x) = g(x) + \dfrac{x_f - x}{x_f - x_0}(y_0 - g_0) + \dfrac{x - x_0}{x_f - x_0}(y_f - g_f),
\end{equation}
which represents all possible functions satisfying the boundary value constraints. Furthermore, the derivatives follow,
\begin{equation}\label{eq:gen_dervs}
    \begin{cases}
    y'(x) &= g'(x) - \dfrac{y_0 - g_0}{x_f - x_0} + \dfrac{y_f - g_f}{x_f - x_0} \\
    y''(x) &= g''(x) \\
    & \; \vdots \\
    y^{(n)}(x) &= g^{(n)}(x) \\
    \end{cases}
\end{equation}
By substituting Eq. (\ref{eq:ex_ce_solved}) and Eq. (\ref{eq:gen_dervs}) into Eq. (\ref{eq:gen_DE}), the DE is transformed to an new DE we define as $\tilde{F}$, which is only a function of the independent variable $x$ and the free-function $g(x)$ where,
\begin{equation*}
    \tilde{F}\left(x,g,g',g''\right) = 0.
\end{equation*}
This DE has a unique solution, it is subject to \emph{no constraints}, and the boundary-value constraints are satisfied regardless of the value of $g(x)$. In order to solve this problem numerically, we express the function $g (x)$ as a linear combination of some known basis function,
\begin{equation}\label{eq:basis}
    g (x) = \B{\xi} \T\B{h} (z) \quad \text{where} \quad z = z(x)
\end{equation}
where $\B{\xi}$ is a $m \times 1$ vector of unknown coefficients (where $m$ is the number of basis functions) and $\B{h}(z)$ is the vector of the selected basis. In general, the basis functions are defined on a specific domain (Chebyshev and Legendre polynomials are defined on $z\in[-1,+1]$, Fourier series is defined on $z\in[-\pi,+\pi]$, etc.) so these functions must be linearly mapped to the independent variable $x$. This can be done using the equations, 
\begin{equation*}
z = z_0 + \dfrac{z_f-z_0}{x_f-x_0}(x - x_0) \qquad \leftrightarrow \qquad x = x_0 + \dfrac{x_f-x_0}{z_f-z_0}(z - z_0).
\end{equation*}
The subsequent derivatives that the free-function defined in Eq. (\ref{eq:basis}) follow,
\begin{equation*}
    \dfrac{\dd^{n} g}{\dd x^{n}} = \B{\xi} \T  \dfrac{\dd^{n} \B{h}(z)}{\dd z^{n}} \left(\dfrac{\dd z}{\dd x}\right)^{n},
\end{equation*}
where by defining,
\begin{equation*}
c := \dfrac{\dd z}{\dd x} = \dfrac{z_f - z_0}{x_f - x_0}
\end{equation*}
the expression can be simplified to, 
\begin{equation}\label{eq:basis_derv}
    \dfrac{\dd^{n} g}{\dd x^{n}} = c^{n} \B{\xi} \T  \dfrac{\dd^{n} \B{h}(z)}{\dd z^{n}},
\end{equation}
which defines all mappings of the free-function. Following these steps, the DE becomes a function of only the unknown $\B{\xi}$ parameters and the independent variable $x$,
\begin{equation*}
    \tilde{F}(x, \B{\xi}) = 0,
\end{equation*}
This function is an \emph{unconstrained} optimization problem (linear or nonlinear) that needs to be solved for the unknown parameters $\B{\xi}$. This is done by discretizing over the domain $x\in[x_0,x_f]$ by $N$ points. In this paper, and prior papers we consider the linear basis $\B{h}(z)$ as Chebysehv or Legendre polynomials. The optimal distribution of points is provided by Chebysehv-Gauss-Lobatto collocation points \cite{Colloc,ChebCol}, defined as,
\begin{equation}\label{eq:collo}
    z_k = -\cos\left(\dfrac{k \pi}{N}\right) \qquad \text{for} \qquad k = 1, 2, \cdots, N.
\end{equation}
As compared to the uniform distribution point, the collocation point distribution are the most efficient distribution when using orthogonal polynomials in least-squares applications. This results in a better condition number as the number of basis functions, $m$, increases. Once discretized, many optimization schemes can be applied, but for this paper scaled QR decomposition is used. 

Slight adjustments must be made in order to apply the TFC method to hybrid systems, but the derivation in the following sections leverages the framework of the example presented above and remains within the TFC framework.

\section{Segment approach to TFC}

The simplest example of a hybrid system is a DE with a discrete jump in the dynamic behavior at a single point along the domain. When solving a two-point BVP according to these dynamics, no only must the solution satisfy the boundary condition, but it must also preserve the $C^1$ continuity over the jump.

\begin{figure}[ht]
    \centering\includegraphics[width=.65\linewidth]{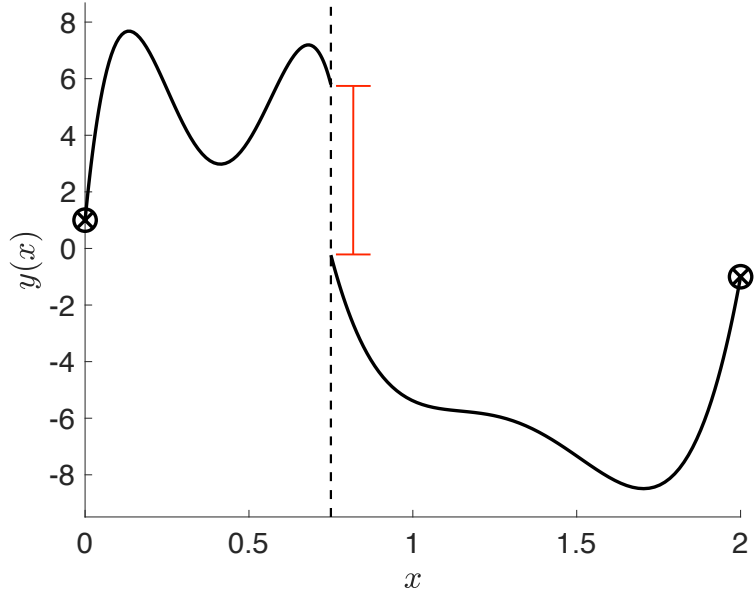}
    \caption{Using a multiple shooting method, the residual of continuity needs to be minimized to solve a hybrid BVP. However, using the TFC approach an analytical expression can be derived that satisfies all constraints of the problem without the need for an outer loop.}
    \label{fig:seg2_notoc}
\end{figure}

Therefore, let us consider a function where it is desired that $C^1$ continuity is also a constraint along with the boundary point constraints. Classical approaches using a shooting method can be used, but in this the residual of the junction must be minimized through an outer optimization loop. A summary of this is shown in Fig. \ref{fig:seg2_notoc}, where the bar relates to the residual of the continuity constraint. In prior studies of TFC \cite{U-ToC,LDE,NDE}, a least-squares solution of BVPs was detailed which removed the need for shooting methods. 
This can be leveraged for the inclusion of continuity constraints and ultimately provide a method to solve BVPs arising from hybrid systems. In order to solve the problem presented in Fig. \ref{fig:seg2_notoc}, we look to develop an analytical expression such that both boundary-value constraints are satisfied along with the function and derivative continuities between the two sub-intervals. For this problem, let $x \in [x_0, x_f]$, and let $x_1$ define the switching of the dynamics at a point along the domain. Additionally, let us indicate with $y_{(1)}(x)$ where $x \in [x_0, x_1]$ and $y_{(2)}(x)$ where $x \in [x_1, x_f]$ the two searched solutions in the two subsequent intervals. In the following section, two distinct approaches are developed and analyzed for ease of implementation and scalability. 

\subsection{Embedded continuity (cascade) method}
In this approach, we develop a method that embeds the junctions as relative constraints shared by each segment\footnote{In this paper, the segment is denoted by a subscript in parentheses to distinguish it from other notations. For example, the function over the first segment is denoted as $y_{(1)}$.}. First, consider the two segments defined by separate \ces\,
\begin{eqnarray}
    y_{(1)}(x) &= g_{(1)}(x) + \eta_1 + \eta_2 x \label{eq:cas1}\\
    y_{(2)}(x) &= g_{(2)}(x) +\eta_3 + \eta_4 x \label{eq:cas2}
\end{eqnarray}
with the constraints, $y (x_0) = y_0$, $y_{(1)} (x_1) = y_{(2)} (x_1) = y_1$, and $y (x_f) = y_f$. Applying these constraints to the two \ces\, leads to
\begin{equation*}
    \begin{aligned}
    y_0 &= g_{(1)}(x_0) + \eta_1 + \eta_2 x_0 \\ y_1 &=  g_{(1)}(x_1) + \eta_1 + \eta_2 x_1 \\ y_1 &= g_{(2)} (x_1) + \eta_3 + \eta_4 x_1 \\
    y_f &= g_{(2)}(x_f) + \eta_3 + \eta_4 x_f\end{aligned}
\end{equation*}
which can be organized into two systems of equations which can be solved for the unknown $\eta_k$ terms.
\begin{equation*}
\begin{aligned}
\begin{Bmatrix} y_0 - g_{(1)}(x_0) \\ y_1 - g_{(1)}(x_1)\end{Bmatrix} &= \begin{bmatrix} 1 & x_0 \\ 1 & x_1\end{bmatrix}\begin{Bmatrix} \eta_1 \\ \eta_2 \end{Bmatrix} \\
\begin{Bmatrix} y_1 - g_{(2)}(x_1) \\ y_f - g_{(2)}(x_f)\end{Bmatrix} &= \begin{bmatrix} 1 & x_1 \\ 1 & x_f\end{bmatrix}\begin{Bmatrix} \eta_3 \\ \eta_4 \end{Bmatrix},
\end{aligned}
\end{equation*}
This system can be solved through matrix inversion leading to the coefficients, 
\begin{equation*}
    \begin{aligned}
    \eta_1 &= \dfrac{1}{x_1-x_0} \Big[x_1\Big(y_0 - g_{(1)}(x_0)\Big) - x_0\Big(y_1 - g_{(1)}(x_1)\Big)\Big]\\
    \eta_2 &= \dfrac{1}{x_1-x_0} \Big[\Big((y_1 - g_{(1)}(x_1)\Big) - \Big(y_0 - g_{(1)}(x_0)\Big)\Big] \\
    \eta_3 &= \dfrac{1}{x_f-x_1} \Big[x_f\Big(y_1 - g_{(2)}(x_1)\Big) - x_1\Big(y_f - g_{(2)}(x_f)\Big)\Big]\\
    \eta_4 &= \dfrac{1}{x_f-x_1} \Big[\Big((y_f - g_{(2)}(x_f)\Big) - \Big(y_1 - g_{(2)}(x_1)\Big)\Big]
    \end{aligned}
\end{equation*}
Substituting these expressions into the original \ces\ defined by Eq. (\ref{eq:cas1}) and Eq. (\ref{eq:cas2}), the functions $y_{(1)}(x)$ and $y_{(2)}(x)$ become,
\begin{equation*}
y_{(1)}(x) = g_{(1)}(x) + \underbrace{\dfrac{x_1 - x}{x_1 - x_0}}_\text{\normalsize$\alpha_1(x_{(1)})$}\Big(y_0 - g_{(1)}(x_0)\Big) + \underbrace{\dfrac{x - x_0}{x_1 - x_0}}_\text{\normalsize$\alpha_2(x_{(1)})$}\Big(y_1 - g_{(1)}(x_1)\Big),
\end{equation*}
and
\begin{equation*}
    y_{(2)}(x) = g_{(2)}(x) + \underbrace{\dfrac{x_f - x}{x_f - x_1}}_\text{\normalsize$\alpha_1(x_{(2)})$}\Big(y_1 - g_{(2)}(x_1)\Big) + \underbrace{\dfrac{x - x_1}{x_f - x_1}}_\text{\normalsize$\alpha_2(x_{(2)})$}\Big(y_f-g_{(2)}(x_f)\Big).
\end{equation*}
In these equations, the subscript in parentheses on the independent variable indicates the segment considered. For example, the notation of $x_{(1)}$ denotes $x \in [x_0, x_1]$ (or $x$ over the first segment). 

An attentive reader will notice that the quantities defined by $\alpha_1(x)$ and $\alpha_2(x)$ are solely functions of the independent variable. These functions act as continuous switching functions which are dependent on the constraint locations. A summary of the initial-final value switching functions and their derivatives is given in Table \ref{tab:alpha}.
\begin{table}[h]
\begin{center}
\begin{tabular}{SSSS} \toprule
    {} & {Initial Value Switching Function} & {Initial Value Switching Function} \\
    {} & {$\alpha_1(x)$} & {$\alpha_2(x)$} \\ \bottomrule \midrule 
    {$(\cdot)$} & {$\dfrac{x_f - x}{\Delta x}$} & {$\dfrac{x - x_0}{\Delta x}$}\\ \midrule
    {$\dfrac{\dd}{\dd x}(\cdot)$} & {$\dfrac{-1}{\Delta x}$} & {$\dfrac{1}{\Delta x}$} \\ \midrule
    {$\dfrac{\dd^2}{\dd x^2}(\cdot)$} & {$0$} & {$0$} \\ \midrule
 \bottomrule
\end{tabular}
\end{center}
\caption{Initial-final value switching functions defined for a general domain of $x \in [x_0, x_f]$. The switching function defined by $\alpha_1(x)$ multiples the initial value constraint and $\alpha_1(x)$ multiples the final value constraint.}
\label{tab:alpha}
\end{table}

Simplifying the expressions we still have the unknown in $y_1$,
\begin{equation}\label{eq:ce1}
y_{(1)}(x) = g_{(1)}(x) + \alpha_1(x_{(1)}) \Big(y_0 - g_{(1)}(x_0)\Big) + \alpha_2(x_{(1)})\Big(y_1 - g_{(1)}(x_1)\Big),
\end{equation}
and
\begin{equation}\label{eq:ce2}
    y_{(2)}(x) = g_{(2)}(x) + \alpha_1(x_{(2)})\Big(y_1 - g_{(2)}(x_1)\Big) + \alpha_2(x_{(2)})\Big(y_f-g_{(2)}(x_f)\Big).
\end{equation}
which can be solved for by taking a derivative of Eq. (\ref{eq:ce1}) and Eq. (\ref{eq:ce2}) and setting them equal (thus enforcing the continuity of the derivative over the jump). This produces a linear function in the unknown,
\begin{equation*}
\begin{aligned}
 g'_{(1)}(x_1) + &\alpha'_1(x_{(1)}(x_1)) \Big(y_0 - g_{(1)}(x_0)\Big) + \alpha'_2(x_{(1)}(x_1))\Big(y_1 - g_{(1)}(x_1)\Big) \\ = &g'_{(2)}(x_1) + \alpha'_1(x_{(2)}(x_1))\Big(y_1 - g_{(2)}(x_1)\Big) + \alpha'_2(x_{(2)}(x_1))\Big(y_f-g_{(2)}(x_f)\Big)
 \end{aligned}
\end{equation*}
which can be solved to obtain
\begin{equation*}
\begin{split}
    y_1 =& \dfrac{g'_{(2)}(x_1) - g'_{(1)}(x_1) - \alpha'_1(x_{(1)}(x_1)) \Big(y_0 - g_{(1)}(x_0)\Big) + \alpha'_2(x_{(1)}(x_1))g_{(1)}(x_1)}{\alpha'_2(x_{(1)}(x_1)) - \alpha'_1(x_{(2)}(x_1))} + \\
    ~&+\dfrac{- \alpha'_1(x_{(2)}(x_1))g_{(2)}(x_1) + \alpha'_2(x_{(2)}(x_1))\Big(y_f - g_{(2)}(x_f)\Big)}{\alpha'_2(x_{(1)}(x_1)) - \alpha'_1(x_{(2)}(x_1))}.
\end{split}
\end{equation*}
The function describing $y_1$ remains linear in the unknown parameters, however, the terms shows up in both $y_{(1)}$ and $y_{(2)}$ and thus each functions shares information from the other segment. If this extended to more segments, a cascade effect will occur where each segment is dependent on multiple other segments. This occurs because this formulation is capturing both the original constraints and continuity constraints using the free functions. To avoid this issue (where the matrix is fully populated), another approach is devised by treating the value of the junctions as another free parameter to be determined. This has the major advantage of scalability to larger systems and sparsity of the system, at the cost of a few more variables to estimate.

\subsection{Unknown relative parameters method}\label{sec:rel_par}
\begin{figure}[ht]
    \centering\includegraphics[width=.65\linewidth]{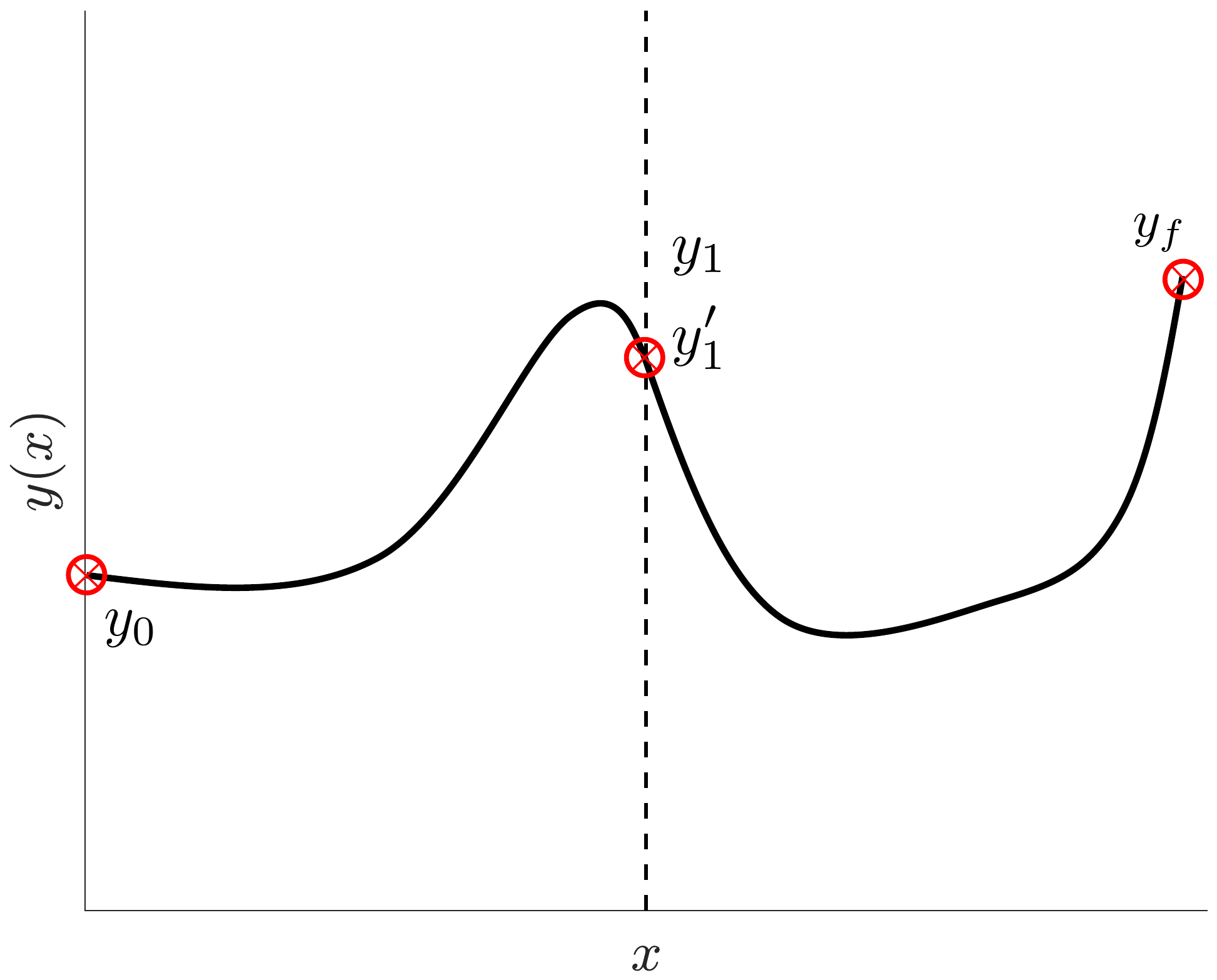}
    \caption{Illustration of piecewise TFC approach enforcing $C^1$ continuity over two segments.}
    \label{fig:seg2_template}
\end{figure}

Again, we look to solve the same problem detailed in Fig. \ref{fig:seg2_notoc}. However, the \ce\ will be different because we will enforce the continuity of the derivative in the \ce. A visualization of this approach is detailed in Fig. \ref{fig:seg2_template}.
The constrained expression for both segments are now searched using,
\begin{equation}\label{eq:ce_rp1}
    y_{(1)}(x) = g_{(1)}(x) + \eta_1 + \eta_2 x + \eta_3 x^2
\end{equation}
\begin{equation}\label{eq:ce_rp2}
    y_{(2)}(x) = g_{(2)}(x) + \eta_4 + \eta_5 x + \eta_6 x^2.
\end{equation}
By applying the constraints, the following conditions are, 
\begin{equation*}
    \begin{aligned}
    y_0 &= g_{(1)}(x_0) + \eta_1 + \eta_2 x_0 + \eta_3 x_0^2 \\
    y_1 &= g_{(1)}(x_1) + \eta_1 + \eta_2 x_1 + \eta_3 x_1^2 \\
    y'_1 &= g'_{(1)}(x_1) + \eta_2 + 2 \eta_3 x_1 \\
    y_1 &= g_{(2)}(x_1) + \eta_4 + \eta_5 x_1 + \eta_6 x_1^2 \\
    y'_1 &= g'_{(2)}(x_1) + \eta_5 + 2 \eta_6 x_1 \\
    y_f &= g_{(2)}(x_f) + \eta_4 + \eta_5 x_f + \eta_6 x_f^2
    \end{aligned}
\end{equation*}
which leads to the two systems of equations
\begin{equation*}
    \begin{aligned}
    \begin{Bmatrix} y_0 - g_{(1)}(x_0) \\  y_1 - g_{(1)}(x_1) \\ y'_1 - g'_{(1)}(x_1) \end{Bmatrix} &= \begin{bmatrix} 1 & x_0 & x_0^2 \\ 1 & x_1 & x_1^2 \\ 0 & 1 & 2x_1\end{bmatrix} \begin{Bmatrix} \eta_1 \\ \eta_2 \\ \eta_3 \end{Bmatrix} \\
    \begin{Bmatrix} y_1 - g_{(2)}(x_1) \\  y'_1 - g'_{(2)}(x_1) \\ y_f - g_{(2)}(x_f) \end{Bmatrix} &= \begin{bmatrix} 1 & x_1 & x_1^2 \\ 0 & 1 & 2x_1 \\ 1 & x_f & x_f^2\end{bmatrix} \begin{Bmatrix} \eta_4 \\ \eta_5 \\ \eta_6 \end{Bmatrix}
    \end{aligned}
\end{equation*}
The unknown $\eta$ values can be solved for by matrix inversion leading to,
\begin{equation*}
    \begin{Bmatrix} \eta_1 \\ \eta_2 \\ \eta_3 \end{Bmatrix} = \dfrac{1}{x_1^2 - 2 x_1 x_0 + x_0^2} \begin{bmatrix} x_1^2 & x_0^2 - 2 x_0 x_1 & x_0 x_1^2 -x_0^2 x_1 \\ -2x_1 & 2 x_1 & x_0^2 - x_1^2\\ 1 & -1 & x_1 - x_0 \end{bmatrix} \begin{Bmatrix} y_0 - g_{(1)}(x_0) \\  y_1 - g_{(1)}(x_1) \\ y'_1 - g'_{(1)}(x_1) \end{Bmatrix}
\end{equation*}
and
\begin{equation*}
    \begin{Bmatrix} \eta_4 \\ \eta_5 \\ \eta_6 \end{Bmatrix} = \dfrac{1}{x_f^2 - 2x_f x_1 + x_1^2} \begin{bmatrix} x_f^2 -2 x_1 x_f & x_1^2 x_f - x_1 x_f^2 & x_1^2 \\ 2x_1 & x_f^2 - x_1^2 & -2x_1 \\ -1 & x_1 - x_f & 1 \end{bmatrix} \begin{Bmatrix} y_1 - g_{(2)}(x_1) \\  y'_1 - g'_{(2)}(x_1) \\ y_f - g_{(2)}(x_f) \end{Bmatrix}
\end{equation*}
Plugging these into Eqs. (\ref{eq:ce_rp1}) and (\ref{eq:ce_rp2}), the \ce\ takes a similar form to the prior section with the $\beta(x)$ function representing the switching functions based on the independent variable similar to the $\alpha(x)$ functions introduced earlier. A summary of these $\beta(x)$ functions and derivatives are provided in Table \ref{tab:beta13} and \ref{tab:beta46}. 
\begin{table}[H]
\begin{center}
\begin{tabular}{SSSSS} \toprule
{} & {Initial Value Switching Function} & {Final Value Switching Function} & {Final Derivative Switching Function} \\
{} & {$\beta_1(x)$} & {$\beta_2(x)$} & {$\beta_3(x)$} \\ \bottomrule \midrule 
{$(\cdot)$} & {$\dfrac{(x-x_f)^2}{\Delta x^2}$} & {$\dfrac{(x_0 - x) (x + x_0 - 2 x_f)}{\Delta x^2}$} & {$\dfrac{(x-x_0)(x-x_f)}{\Delta x}$} \\ \midrule
{$\dfrac{\dd}{\dd x}(\cdot)$} & {$\dfrac{2(x-x_f)}{\Delta x^2}$} & {$-\dfrac{2(x-x_f)}{\Delta x^2}$} & {$\dfrac{2x - x_0 - x_f}{\Delta x}$} \\ \midrule
{$\dfrac{\dd^2}{\dd x^2}(\cdot)$} & {$\dfrac{2}{\Delta x^2}$} & {$-\dfrac{2}{\Delta x^2}$} & {$\dfrac{2}{\Delta x}$} \\ \midrule
\bottomrule
\end{tabular}
\end{center}
\caption{Switching functions for first segment defined for a general domain of $x \in [x_0, x_f]$.}
\label{tab:beta13}
\end{table}

\begin{table}[H]
\begin{center}
\begin{tabular}{SSSSS} \toprule
{} & {Initial Value Switching Function} & {Initial Derivative Switching Function} & {Final Value Switching Function} \\
{} & {$\beta_4(x)$} & {$\beta_5(x)$} & {$\beta_6(x)$} \\  \bottomrule \midrule
{$(\cdot)$} & {$\dfrac{(x_f-x)(x - 2x_0 + x_f)}{\Delta x^2}$} & {$\dfrac{(x-x_0)(x_f-x)}{\Delta x}$} & {$\dfrac{(x-x_0)^2}{\Delta x^2} $} \\ \midrule
{$\dfrac{\dd}{\dd x}(\cdot)$} & {$\dfrac{-2(x-x_0)}{\Delta x^2}$} & {$\dfrac{-2x+x_0 + x_f}{\Delta x}$} & {$\dfrac{2(x-x_0)}{\Delta x^2}$} \\ \midrule
{$\dfrac{\dd^2}{\dd x^2}(\cdot)$} & {$\dfrac{-2}{\Delta x^2}$} & {$\dfrac{-2}{\Delta x}$} & {$\dfrac{2}{\Delta x^2}$} \\ \midrule
\bottomrule
\end{tabular}
\end{center}
\caption{Switching functions for last segment defined for a general domain of $x \in [x_0, x_f]$.}
\label{tab:beta46}
\end{table}
The constrained expressions become,
\begin{equation}\label{eq:ce_rp1_full}
    y_{(1)}(x) = g_{(1)}(x) + \beta_1(x)\Big(y_0 - g_{(1)}(x_0)\Big) + \beta_2(x)\Big(y_1 - g_{(1)}(x_1)\Big) + \beta_3(x)\Big(y'_1 - g'_{(1)}(x_1)\Big)
\end{equation}
\begin{equation}\label{eq:ce_rp2_full}
    y_{(2)}(x) = g_{(2)}(x) + \beta_4(x)\Big(y_1 - g_{(2)}(x_1)\Big) + \beta_5(x)\Big(y'_1 - g'_{(2)}(x_1)\Big) + \beta_6(x)\Big(y_f - g_{(2)}(x_f)\Big)
\end{equation}
Now, by allowing $y_1$ and $y'_1$ to be free parameters, the two equations defined by Eq. (\ref{eq:ce_rp1_full}) and Eq. (\ref{eq:ce_rp2_full}) are linear in the free parameters. This is the main outcome of the original development of the TFC method in Ref. \cite{U-ToC}. Then, by expressing $g_{(1)}(x)$ and $g_{(2)}(x)$ using Eq. (\ref{eq:basis}), these \ces\ can be used in the solution of a differential equation with a jump in dynamics as $x_1$. The benefit of this method is it can be generalized for ``$n$'' jumps in dynamics which leads to a block diagonal system of equations. The generalization is provided in the following section.
\section{Generalization for ``$n$'' segments}
\begin{figure}[ht]
    \centering\includegraphics[width=.65\linewidth]{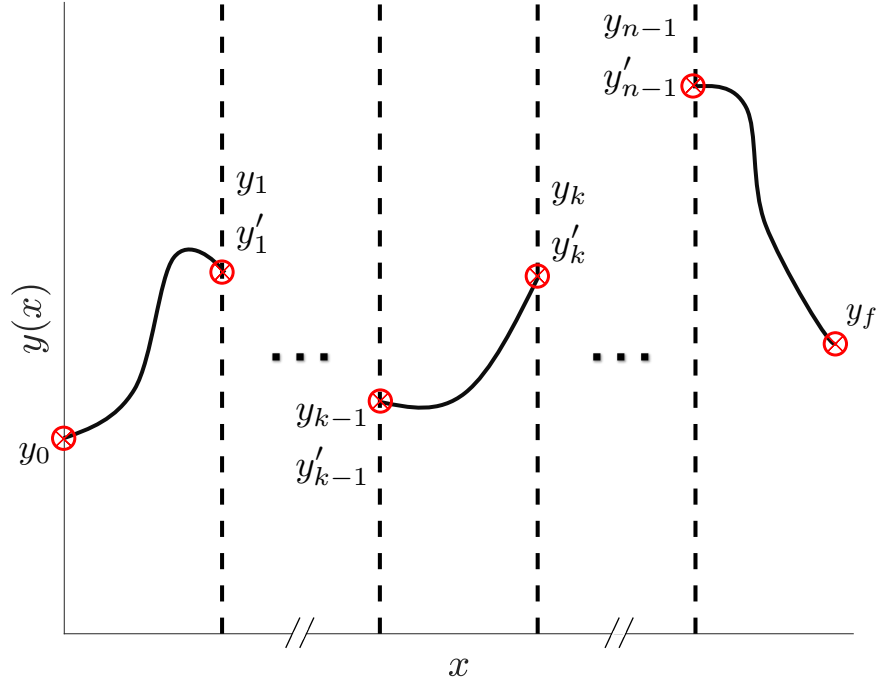}
    \caption{Illustration of segment TFC approach to enforcing $C^1$ continuity over ``$n$'' segments.}
    \label{fig:n_seg}
\end{figure}

Suppose the problem is subject to ``$n$'' jumps in dynamics as detailed in Fig. \ref{fig:n_seg}. This case is the generalization of the problem presented in Section \ref{sec:rel_par}. 
Additionally, this generalization necessitates the introduction of another set of switching functions that can be derived using the TFC method. By analyzing Fig. \ref{fig:n_seg}, it can be seen that in the interior segments (the segments only governed by the constraints on continuity) the number of constraint terms becomes four. The \ce\ for this constraint type produces the \ce\,
\begin{equation*}
    y(x) = g(x) + \gamma_1(x) \Big(c_0 - g_0\Big) + \gamma_2(x) \Big(c_f - g_f\Big) + \gamma_3(x) \Big(c'_0 - g'_0\Big) + \gamma_4(x) \Big(c'_f - g'_f\Big)
\end{equation*}
where $c_0$, $c_f$, $c'_0$, and $c'_f$ are the values and derivative constraints respectively, and the $\gamma(x)$ functions are similar to the switching functions of $\alpha(x)$ and $\beta(x)$. The values of the $\gamma(x)$ switching functions and their derivatives are provided in Tables \ref{tab:gam12} and \ref{tab:gam34}.

\begin{table}[H]
\begin{center}
\begin{tabular}{SSS} \toprule
{} & {Initial Value Switching Function} & {Final Value Switching Function}\\
{} & {$\gamma_1(x)$} & {$\gamma_2(x)$} \\ \bottomrule \midrule
{$(\cdot)$} & {$1 + \dfrac{2(x-x_0)^3}{\Delta x^3} - \dfrac{3 (x-x_0)^2}{\Delta x^2}$} & {$- \dfrac{2 (x-x_0)^3}{\Delta x^3} + \dfrac{3 (x-x_0)^2}{\Delta x^2}$} \\ \midrule
{$\dfrac{\dd}{\dd x}(\cdot)$} & {$\dfrac{6(x-x_0)^2}{\Delta x^3} - \dfrac{6(x-x_0)}{\Delta x^2}$} & {$-\dfrac{6(x-x_0)^2}{\Delta x^3} + \dfrac{6(x-x_0)}{\Delta x^2}$} \\ \midrule
{$\dfrac{\dd^2}{\dd x^2}(\cdot)$} & {$\dfrac{12(x-x_0)}{\Delta x^3} - \dfrac{6}{\Delta x^2}$} & {$-\dfrac{12(x-x_0)}{\Delta x^3} + \dfrac{6}{\Delta x^2}$}\\ \midrule
\bottomrule
\end{tabular}
\end{center}
\caption{Switching functions for function constraints for middle segments defined for a general domain of $x \in [x_0, x_f]$.}
\label{tab:gam12}
\end{table}

\begin{table}[H]
\begin{center}
\begin{tabular}{SSS} \toprule
{} & {Initial Derivative Switching Function} & {Final Derivative Switching Function}\\
{}& {$\gamma_3(x)$} & {$\gamma_4(x)$} \\ \bottomrule \midrule
{$(\cdot)$} & {$(x-x_0) + \dfrac{(x-x_0)^3}{\Delta x^2} - \dfrac{2 (x-x_0)^2}{\Delta x}$} & {$\dfrac{(x-x_0)^3}{\Delta x^2} - \dfrac{(x-x_0)^2}{\Delta x}$} \\ \midrule
{$\dfrac{\dd}{\dd x}(\cdot)$} & {$1 + \dfrac{3(x-x_0)^2}{\Delta x^2} - \dfrac{4(x-x_0)}{\Delta x}$} & {$\dfrac{3(x-x_0)^2}{\Delta x^2} - \dfrac{2(x-x_0)}{\Delta x}$} \\ \midrule
{$\dfrac{\dd^2}{\dd x^2}(\cdot)$} & {$\dfrac{6(x-x_0)}{\Delta x^2} - \dfrac{4}{\Delta x}$} & {$\dfrac{6(x-x_0)}{\Delta x^2} - \dfrac{2}{\Delta x}$}\\ \midrule
\bottomrule
\end{tabular}
\end{center}
\caption{Switching functions for derivative constraints for middle segments defined for a general domain of $x \in [x_0, x_f]$.}
\label{tab:gam34}
\end{table}
Therefore, the \ces\ for ``$n$'' segments can be generally written by three unique equations. The first segment is identical to Eq. (\ref{eq:ce_rp1_full}),
\begin{equation*}
    y_{(1)}(x) = g_{(1)}(x) + \beta_1(x)\Big(y_0 - g_{(1)}(x_0)\Big) + \beta_2(x)\Big(y_1 - g_{(1)}(x_1)\Big) + \beta_3(x)\Big(y'_1 - g'_{(1)}(x_1)\Big),
\end{equation*}
The middle segments utilize the switching functions developed in this section since these segments have relative constraints on both ends,
\begin{equation*}
\begin{aligned}
    y_{(k)}(x) = g_{(k)}(x) + \gamma_1(x) \Big(y_{k-1} - g_{(k)}(x_{k-1}) \Big) &+ \gamma_2(x) \Big(y_{k} - g_{(k)}(x_{k}) \Big) \\ &+ \gamma_3(x) \Big(y'_{k-1} - g'_{(k)}(x_{k-1}) \Big) + \gamma_4(x) \Big(y'_{k} - g'_{(k)}(x_{k}) \Big),
\end{aligned}
\end{equation*}
where $k = 2, 3, \cdots, n-1$. The final segment is similar to Eq. (\ref{eq:ce_rp2_full}), except the initial conditions of the final segment are based on the final conditions of segment $n-1$.
\begin{equation*}
    y_{(n)}(x) = g_{(n)}(x) + \beta_4(x)\Big(y_{n-1} - g_{(n)}(x_{n-1})\Big) + \beta_5(x)\Big(y'_{n-1} - g'_{(n)}(x_{n-1})\Big) + \beta_6(x)\Big(y_f - g_{(n)}(x_f)\Big).
\end{equation*}
By expressing the free functions $g(x)$ according to Eq. (\ref{eq:basis}) and  Eq. (\ref{eq:basis_derv}) and discretizing the domains, the generalization can be written in a compact block diagonal matrix of the form,
\begin{equation*}
    \B{y} = \mathbb{A} \, \Xi + \mathbb{B}
\end{equation*}
where
\begin{equation*}
    \mathbb{A} = \begin{bmatrix} R_1 \\ \vdots \\ R_k \\ \vdots \\ R_n \end{bmatrix} \quad \text{such that:} \quad \begin{cases} \, R_1 &= \, \begin{bmatrix} H_{(1)} & \beta_2(\B{x}_{(1)}) & \beta_3(\B{x}_{(1)}) & \B{0}_{N \times [m+(m+2)(n-2)]} \end{bmatrix} \\
    \, R_k &= \, \begin{bmatrix} \B{0}_{N \times [m + (m+2)(k-2)]} & \gamma_1(\B{x}_{(k)}) & \gamma_2(\B{x}_{(k)}) & H_{(k)} & \gamma_3(\B{x}_{(k)}) & \gamma_4(\B{x}_{(k)}) & \B{0}_{N \times [m + (m+2)(n-k-1)]} \end{bmatrix} \\
    \, R_n &=  \, \begin{bmatrix} \B{0}_{N \times [m+(m+2)(n-2)]} & \beta_4(\B{x}_{(n)}) & \beta_5(\B{x}_{(n)}) & H_{(n)} \end{bmatrix} \end{cases}
\end{equation*}
where $\Xi$ is a vector of all the free parameters,
\begin{equation*}
    \Xi = \begin{Bmatrix} \B{\xi}_{(1)}\T & y_1 & y'_1 & \hdots &  y_{k-1} & y'_{k-1} & \B{\xi}_{(k)}\T & y_k & y'_k & \hdots & y_{n-1} & y'_{n-1} & \B{\xi}_{(n)}\T  \end{Bmatrix}\T 
\end{equation*}
and $\mathbb{B}$ is a vector with the boundary constraints,
\begin{equation*}
    \mathbb{B} = \begin{Bmatrix}\beta_1(\B{x}_{(1)})y_0 \\ \B{0}_{(N-2) \times 1} \\ \beta_6(\B{x}_{(n)})y_f \end{Bmatrix}.
\end{equation*}
For these equations, $\B{x}_{(k)}$ denotes a discretized vector of x-values such that $\B{x}_{(k)} \in [x_{k-1}, x_k]$ and $H_{(k)}$ is a matrix of the terms multiplied by $\B{\xi}_{(k)}$. For example, the term $H_{(1)}$ is simply,
\begin{equation*}
    H_{(1)} = \Big( \B{h}_{(1)}(\B{x}\T_{(1)}) - \B{h}_{(1)}(x_0)\beta_1(\B{x}\T_{(1)}) -\B{h}_{(1)}(x_1)\beta_2(\B{x}\T_{(1)})  - \B{h}'_{(1)}(x_1)\beta_3(\B{x}\T_{(1)}) \Big)\T
\end{equation*}
Since this is a linear set of equations all subsequent derivatives are the derivatives of the individual components. Making sure that the mapping of the basis is accounted for with Eq. (\ref{eq:basis_derv}), the derivatives order $d$ of $\B{y}$ becomes, 
\begin{equation*}
    \B{y}^{(d)} = \mathbb{A}^{(d)} \, \Xi + \mathbb{B}^{(d)},
\end{equation*}
which is also block diagonal.

\section{Application to Hybrid System BVPs}
The following sections apply the theory developed in the previous sections to three hybrid system scenarios. The first problem considers a hybrid system where the dynamics of the two subsequent segments are both linear DEs. The second scenario considers the transition from a linear system to a nonlinear nonlinear system. The last scenario is dedicated to the nonlinear/nonlinear transition. The proposed approach solves all hybrid system BVPs at machine level accuracy.

\subsection{Linear-Linear Differential Equation Sequence}\label{ex1}

Consider a second-order linear DE  sequence such that,
\begin{equation}\label{eq:linear_linear}
    y''(x) = x^2 + a \quad \text{subject to:} \; \begin{cases} y(0) = 0 \\ y(1) = 1 \end{cases} \text{where} \; \begin{cases} a = 0 \quad & \text{for } x \le 0.5 \\ a = 1 \quad & \text{for } x > 0.5 \end{cases}
\end{equation}
which has a forcing term that is discontinuous at $x = 0.5$. This problem has the analytical solution,
\begin{equation*}
    \left\{\begin{array}{ll} y_{(1)} (x) = \dfrac{1}{12} x^4 + \dfrac{19}{24} x & \text{for} \; x \le 0.5 \\ [8pt] y_{(2)} (x) = \dfrac{1}{12} x^4 + \dfrac{1}{2} x^2 + \dfrac{7}{24} x + \dfrac{1}{8} & \text{for} \; x > 0.5 \end{array}\right.
\end{equation*}
Following the TFC procedure (with additional unknowns) given in the previous sections, this problem can be solved using two segments with continuity enforced at the jump. Specifically, the \ces\ for the two segments are,
\begin{equation*}
    y_{(1)}(x) = g_{(1)}(x) + \beta_1(x)\Big(y_0 - g_{(1)}(x_0)\Big) + \beta_2(x)\Big(y_1 - g_{(1)}(x_1)\Big) + \beta_3(x)\Big(y'_1 - g'_{(1)}(x_1)\Big)
\end{equation*}
\begin{equation*}
    y_{(2)}(x) = g_{(2)}(x) + \beta_4(x)\Big(y_1 - g_{(2)}(x_1)\Big) + \beta_5(x)\Big(y'_1 - g'_{(2)}(x_1)\Big) + \beta_6(x)\Big(y_f - g_{(2)}(x_f)\Big)
\end{equation*}
where $y_1$ and $y'_1$ are the additional unknowns. By defining $g_{(1)}(x)$ and $g_{(2)}(x)$ using Eq. (\ref{eq:basis}), the two \ces\ become,
\begin{equation}\label{eq:ex1_ce1}
    y_{(1)}(x) = \Big( \B{h}_{(1)}(x) - \beta_1(x)\B{h}_{(1)}(x_0) -\beta_2(x) \B{h}_{(1)}(x_1) - \beta_3(x)\B{h}'_{(1)}(x_1) \Big)\T \B{\xi}_{(1)} + \beta_1(x) y_0 + \beta_2(x) y_1 + \beta_3(x) y'_1
\end{equation}
\begin{equation}\label{eq:ex1_ce2}
    y_{(2)}(x) = \Big(\B{h}_{(2)}(x) - \beta_4(x)\B{h}_{(2)}(x_1) -\beta_5(x) \B{h}'_{(2)}(x_1) - \beta_6(x)\B{h}_{(2)}(x_f) \Big) \T \B{\xi}_{(2)} + \beta_4(x) y_1 + \beta_5(x) y'_1 + \beta_6(x) y_f
\end{equation}
whose derivatives are simply
\begin{align}
y'_{(1)}(x) &= \Big( \B{h}'_{(1)}(x) - \beta'_1(x)\B{h}_{(1)}(x_0) -\beta'_2(x) \B{h}_{(1)}(x_1) - \beta'_3(x)\B{h}'_{(1)}(x_1) \Big)\T \B{\xi}_{(1)} + \beta'_1(x) y_0 + \beta'_2(x) y_1 + \beta'_3(x) y'_1 \label{eq:ex1_ce1d} \\
y'_{(2)}(x) &= \Big( \B{h}'_{(2)}(x) - \beta'_4(x)\B{h}_{(2)}(x_1) -\beta'_5(x) \B{h}'_{(2)}(x_1) - \beta'_6(x)\B{h}_{(2)}(x_f) \Big) \T \B{\xi}_{(2)} + \beta'_4(x) y_1 + \beta'_5(x) y'_1 + \beta'_6(x) y_f \label{eq:ex1_ce2d} \\
y''_{(1)}(x) &= \Big( \B{h}''_{(1)}(x) - \beta''_1(x)\B{h}_{(1)}(x_0) -\beta''_2(x) \B{h}_{(1)}(x_1) - \beta''_3(x)\B{h}'_{(1)}(x_1) \Big)\T \B{\xi}_{(1)} + \beta''_1(x) y_0 + \beta''_2(x) y_1 + \beta''_3(x) y'_1 \label{eq:ex1_ce1dd} \\
y''_{(2)}(x) &= \Big(\B{h}''_{(2)}(x) - \beta''_4(x)\B{h}_{(2)}(x_1) -\beta''_5(x) \B{h}'_{(2)}(x_1) - \beta''_6(x)\B{h}_{(2)}(x_f) \Big) \T \B{\xi}_{(2)} + \beta''_4(x) y_1 + \beta''_5(x) y'_1 + \beta''_6(x) y_f \label{eq:ex1_ce2dd}
\end{align}
By plugging Eq. (\ref{eq:ex1_ce1dd}) and Eq. (\ref{eq:ex1_ce2dd}) into the DE given by Eq. (\ref{eq:linear_linear}) and discretizing using the method proved by Eq. (\ref{eq:collo}), a linear system of equations can be formed in terms of the unknown vectors, $\B{\xi}_{(1)}$ and $\B{\xi}_{(2)}$, and the unknown scalars, $y_1$ and $y'_1$. The system is of the form,
\begin{equation}\label{eq:ex1_system}
    \begin{bmatrix} H_{(1)} & \beta''_2(\B{x}_{(1)}) & \beta''_3(\B{x}_{(1)}) & \B{0}_{N \times m} \\ \B{0}_{N \times m} & \beta''_4(\B{x}_{(2)}) & \beta''_5(\B{x}_{(2)}) & H_{(2)} \end{bmatrix}\begin{Bmatrix} \B{\xi}_{(1)} \\ y_1 \\ y'_1 \\ \B{\xi}_{(2)}\end{Bmatrix} = \begin{Bmatrix} \B{x}^2_{(1)} - \beta''_1(\B{x}_{(1)}) y_0 \\ \B{x}^2_{(2)} + \B{1}_{N \times 1} - \beta''_6(\B{x}_{(2)}) y_f \end{Bmatrix}
\end{equation}
where $\B{x}_{(1)}$ is the vector of $x$-values discretized over the first segment $\in [0, 0.5]$, and $\B{x}_{(2)}$ is the vector of $x$-values discretized over the second segment $\in (0.5, 1]$. Both $H_{(1)}$ and $H_{(2)}$ represent a $N \times m$ matrix produced by discretizing the terms multiplied by $\B{\xi}_{(1)}$ and $\B{\xi}_{(2)}$. The system of equations given by Eq. (\ref{eq:ex1_system}) can be solved using any optimization technique. The most straightforward is a least-squares approach used in Refs. \cite{LDE,NDE}. In all examples presented in this study, least-squares is performed using scaled QR decomposition. 

\begin{figure}[h]%
\centering
\subfigure[Solution for the linear-linear DEs sequence.]{%
\label{fig:ex1:function}%
\includegraphics[width=0.475\linewidth]{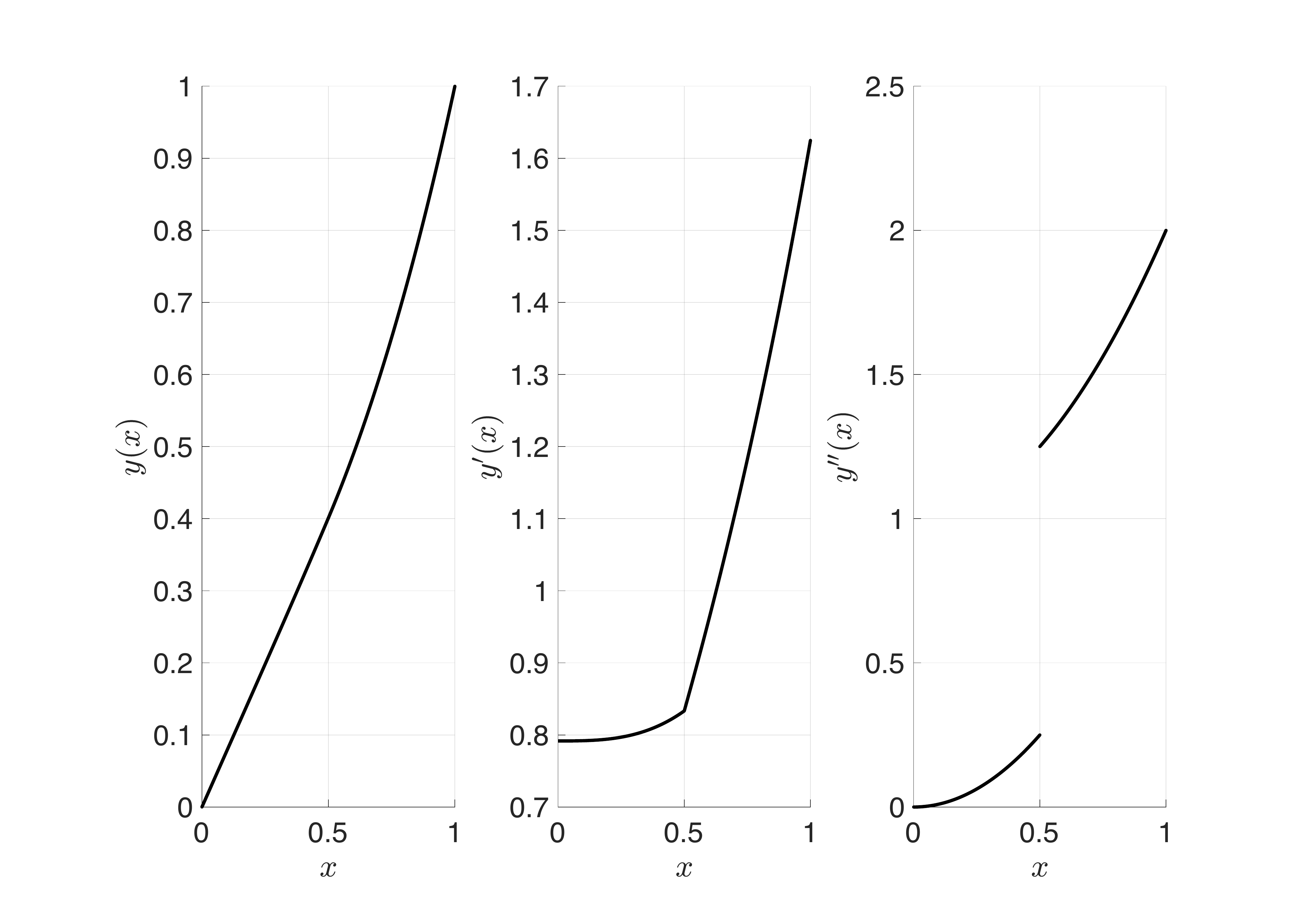}}%
\qquad
\subfigure[Solution error for the linear-linear DEs sequence.]{%
\label{fig:ex1:error}%
\includegraphics[width=0.475\linewidth]{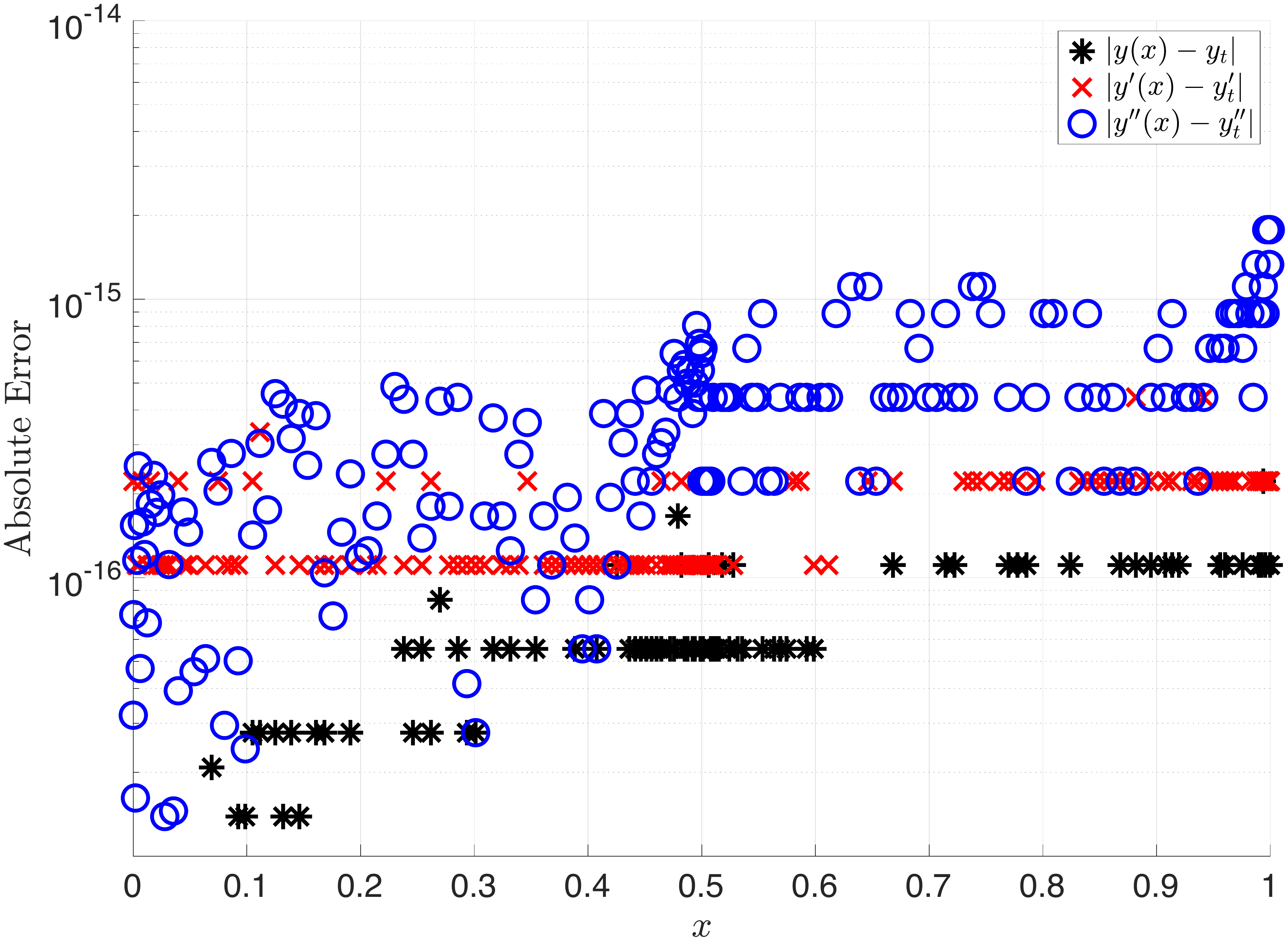}}%
\caption{Results of linear-linear DEs sequence scenario using the TFC approach.}
\label{fig:ex1}
\end{figure}
 By using $N = 100$ points and $m = 3$ basis functions\footnote{In all numerical examples, Chebyshev orthogonal polynomials were used in combination with Legendre-Gauss-Lobatto collocation point distribution.} per segment, this hybrid system was solved using the process described above. Figure \ref{fig:ex1:function} shows the solution of the function and its derivatives while Fig. \ref{fig:ex1:error} displays the absolute error compared to the known analytical solution. It can be seen that the errors are mainly $\mathcal{O}(10^{-15})$. In this specific case, the machine error solution was obtained using just 3 basis functions because the true solution and the basis functions selected were both polynomials. In this case, using more basis function will not not change the results accuracy as all the additional $\xi_k$ terms are estimated all zeros. In the next hybrid systems considered the true solutions are not polynomials.

\subsection{Linear-Nonlinear Differential Equation Sequence}\label{sec:ex2}

\begin{figure}[ht]
    \centering\includegraphics[width=.65\linewidth]{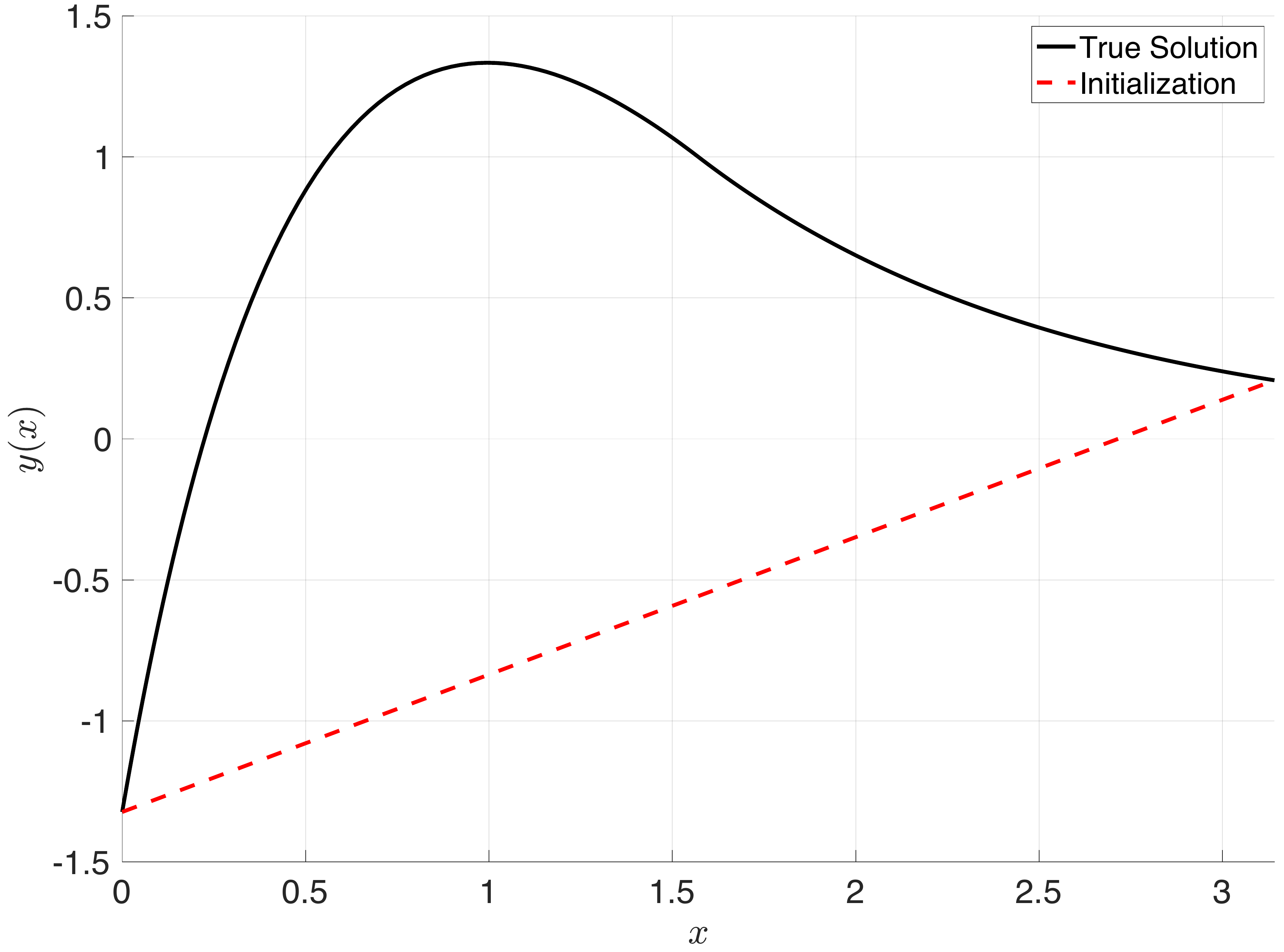}
    \caption{Initial guess and true solution for the linear/nonlinear sequence.}
    \label{fig:ex2_init}
\end{figure}

\begin{figure}[h]%
\centering
\subfigure[Solution of linear-nonlinear differential equation sequence.]{%
\label{fig:ex2:function}%
\includegraphics[width=0.475\linewidth]{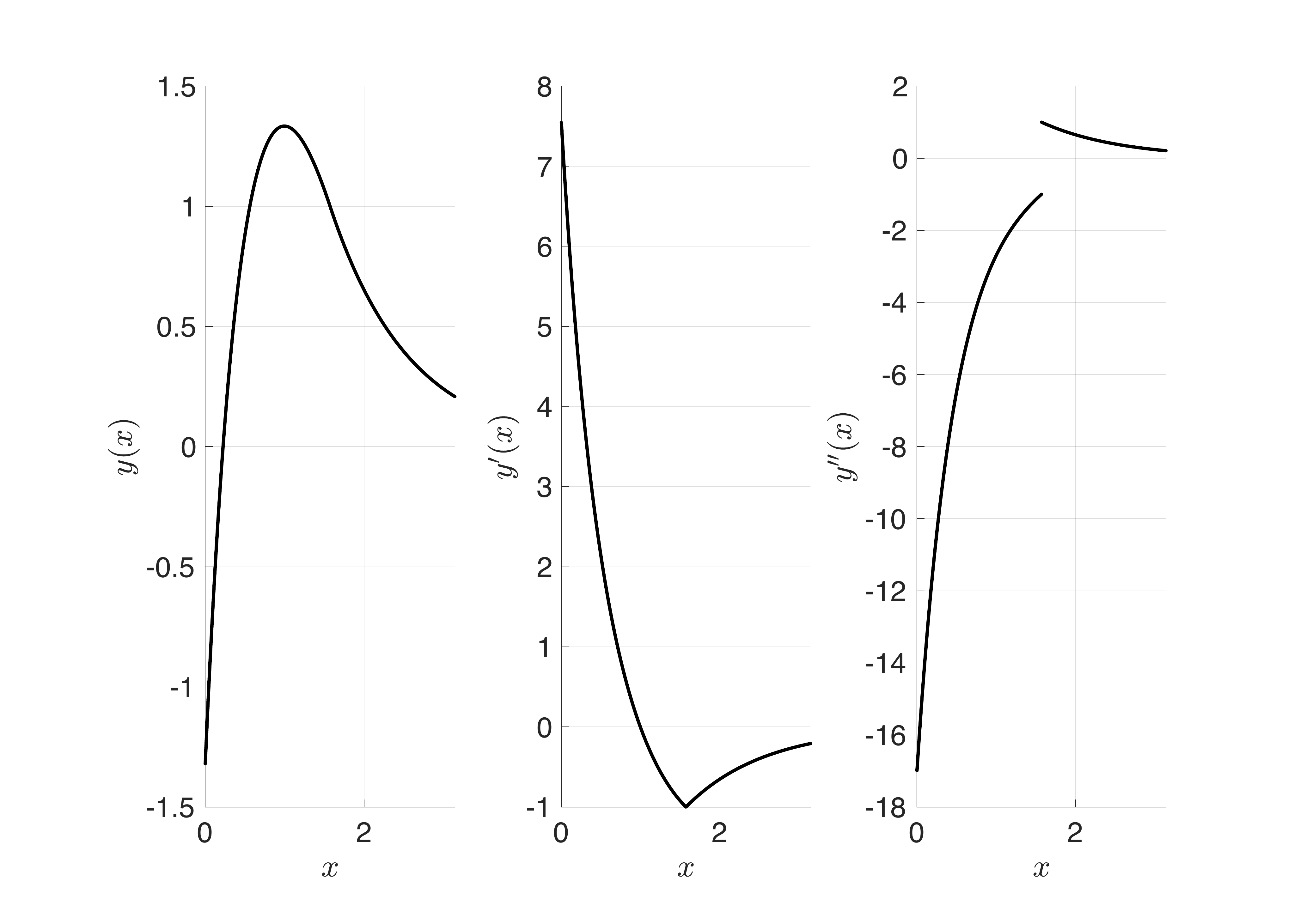}}%
\qquad
\subfigure[Absolute error of solution of linear-nonlinear differential equation sequence.]{%
\label{fig:ex2:error}%
\includegraphics[width=0.475\linewidth]{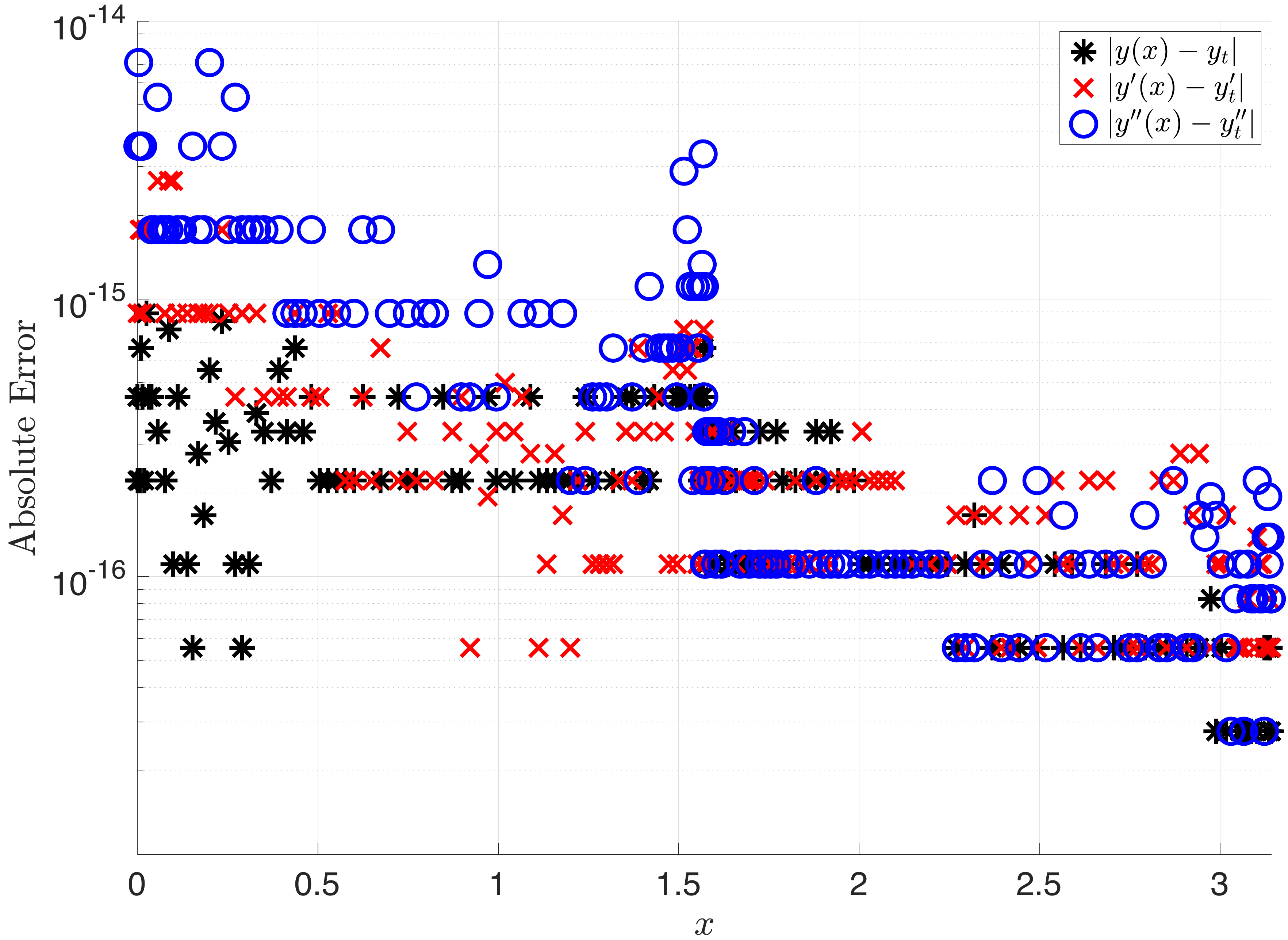}}%
\caption{Results of linear-nonlinear differential equation using the TFC approach. It can be seen that the absolute error for the function and subsequent derivatives are around $10^{-15}$ to $10^{-16}$.}
\label{fig:ex2}
\end{figure}

Consider a second-order linear-nonlinear DE sequence such that,
\begin{equation}\label{eq:linear_nonlinear}
    y'' + y y'^a = e^{\pi/2} - e^{\pi/2 - x} \quad \text{subject to: } \begin{cases} y(0) = \dfrac{9}{10} + \dfrac{1}{10} e^{\pi/2}(5 - 2 e^{\pi/2}) \\ y(\pi) = e^{-\pi/2} \end{cases} \text{and } \begin{cases} a = 0 \quad &\text{for } x \le \pi/2 \\ a = 1 \quad &\text{for } x > \pi/2 \end{cases}
\end{equation}
At the switch, $x = \pi/2$, the the differential equation changes from an linear differential equation to a nonlinear differential equation. This differential equation has the unique solution defined by,
\begin{equation*}
    y(x) = \begin{cases} = - \dfrac{1}{5} e^{\pi - 2 x} + \dfrac{1}{2} e^{\pi/2 - x} + \dfrac{9\cos(x) + 7\sin(x)}{10} \quad &\text{for } x \le \pi/2 \\=  e^{\pi/2 - x} \quad &\text{for } x > \pi/2 \end{cases}
\end{equation*}
The solution of this DE sequence using the TFC approach is similar to the problem solved in Section \ref{ex1}, however, since the sequence has a nonlinear differential equation (over the second segment), an iterative least-squares approach is necessary. For this, we define the residual of the differential equation as the loss function such that,
\begin{equation}\label{eq:L1ex2}
    \mathcal{L}_{(1)}(x) = y''_{(1)}(x) + y_{(1)}(x) - e^{\pi/2} + e^{\pi/2 - x}
\end{equation}
\begin{equation}\label{eq:L2ex2}
    \mathcal{L}_{(2)}(x) = y''_{(2)}(x) + y_{(2)}(x)y'_{(2)}(x) - e^{\pi/2} + e^{\pi/2 - x}
\end{equation}
where $y_{(1)}$, $y''_{(1)}$, $y_{(2)}$, $y'_{(2)}$, and $y''_{(2)}$ are defined by the \ces\ given by Eqs. (\ref{eq:ex1_ce1}-\ref{eq:ex1_ce2dd}) which have the unknown parameters $\B{\xi}_{(1)}$, $y_1$, $y'_1$, and $\B{\xi}_{(2)}$. By substituting these equations into Eqs. (\ref{eq:L1ex2}) and (\ref{eq:L2ex2}) and taking the partials with respect to the unknown parameters a Jacobian can be defined such that, 
\begin{equation}\label{eq:Jacobian}
\mathcal{J} = \begin{bmatrix} \dfrac{\partial \mathcal{L}_{(1)}}{\partial \B{\xi}_{(1)}} & \dfrac{\partial \mathcal{L}_{(1)}}{\partial y_1} & \dfrac{\partial \mathcal{L}_{(1)}}{\partial y'_1} & \B{0}_{N \times m} \\ \B{0}_{N \times m} & \dfrac{\partial \mathcal{L}_{(2)}}{\partial y_1} & \dfrac{\partial \mathcal{L}_{(2)}}{\partial y'_1} & \dfrac{\partial \mathcal{L}_{(2)}}{\partial \B{\xi}_{(2)}} \end{bmatrix}.
\end{equation}
For this problem all terms of Eq. (\ref{eq:Jacobian}) are provided in \ref{sec:linear_nonlinear} in Eqs. (\ref{eq:b1}-\ref{eq:b6}). Defining $\Xi$ as the vector of unknown coefficients such that,
\begin{equation*}
    \Xi = \begin{Bmatrix} \B{\xi}\T_{(1)} & y_1 & y'_1 & \B{\xi}\T_{(2)} \end{Bmatrix}\T,
\end{equation*}
the system can be solved using an iterative least-squares approach where the update to the unknown parameters is given by Newton's method,
\begin{equation*}
    \Xi_{k+1} = \Xi_k - \dd \Xi_k \quad \text{where} \quad \dd \Xi_k = \text{qr}(\mathcal{J}_k, \mathcal{L}_k)
\end{equation*}
Yet, an initial guess must be provided for the iterative least-squares. In the case of BVPs using the TFC method the initial parameters can be determined by connecting the boundary constraints using a straight line. The line initial guess is adopted here in all the numerical BVP tests provided. Therefore, the initial estimate of $y_1$ and $y'_1$ is automatically determined by this initialization. For this problem,
\begin{equation*}
    \Xi_0 = \begin{Bmatrix} \B{0}\T & 1 & -1 & \B{0}\T \end{Bmatrix}\T.
\end{equation*}
A visualization of this initial guess compared to the true solution is provided in Fig. (\ref{fig:ex2_init}).
For the DE presented in Eq. (\ref{eq:linear_nonlinear}), $N = 100$ and
$m = 16$ basis functions were used for each segment. Additionally, the convergence criteria of the nonlinear least-squares method was $|L_2(\mathcal{L})| < \varepsilon$ where $\varepsilon = 10^{-15}$. The solution reached machine error accuracy in $15$ iterations. The results of this numerical test are shown in Figs. \ref{fig:ex2:function} and \ref{fig:ex2:error}. The results show the function, its first two derivatives, and the associated absolute errors compared to the analytical solution. 

\subsection{Nonlinear-Nonlinear Differential Equation Sequence}
\begin{figure}[h]%
\centering
\subfigure[Solution of nonlinear-nonlinear differential equation sequence.]{%
\label{fig:ex3:function}%
\includegraphics[width=0.475\linewidth]{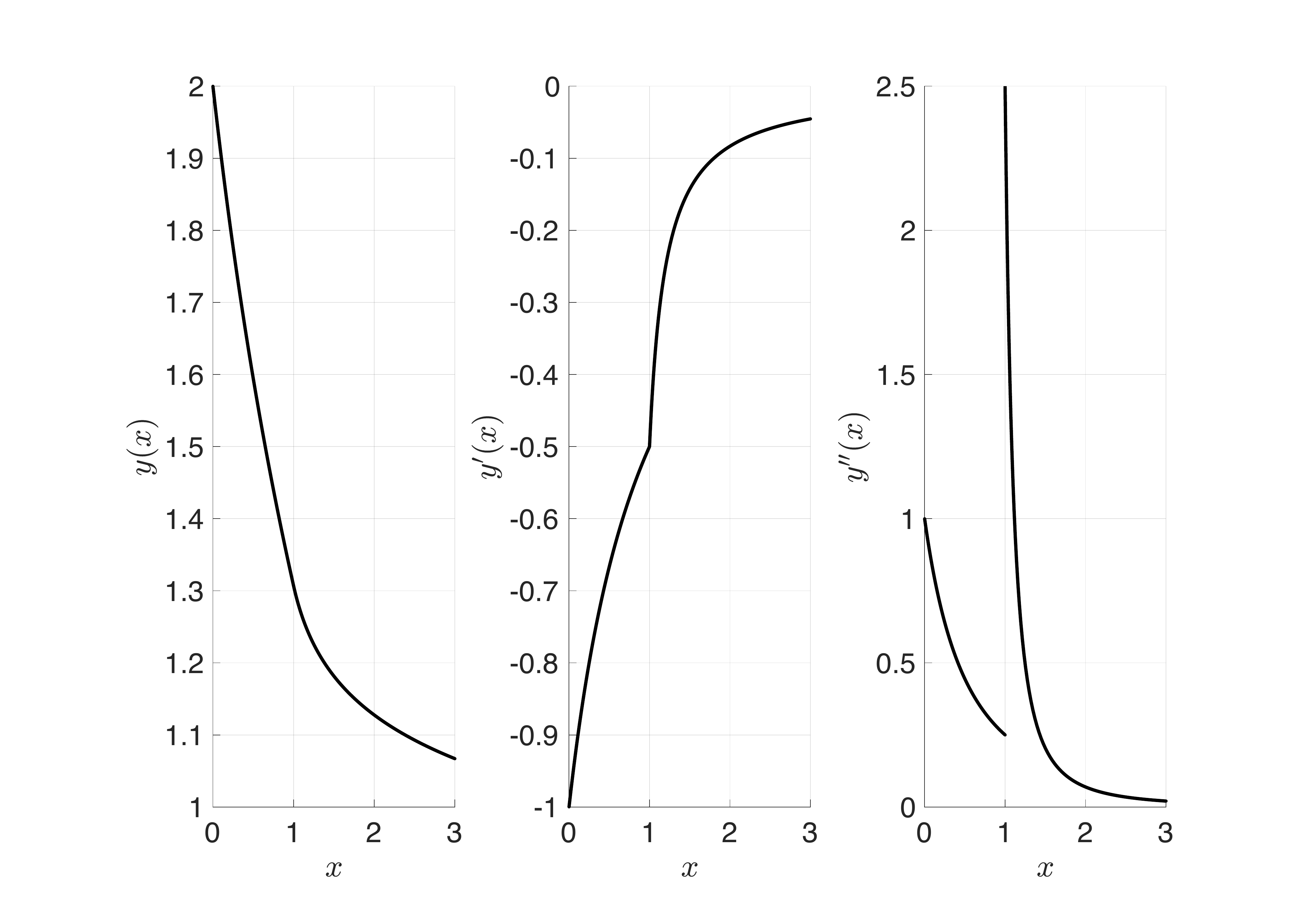}}%
\qquad
\subfigure[Absolute error of solution of nonlinear-nonlinear differential equation sequence.]{%
\label{fig:ex3:error}%
\includegraphics[width=0.475\linewidth]{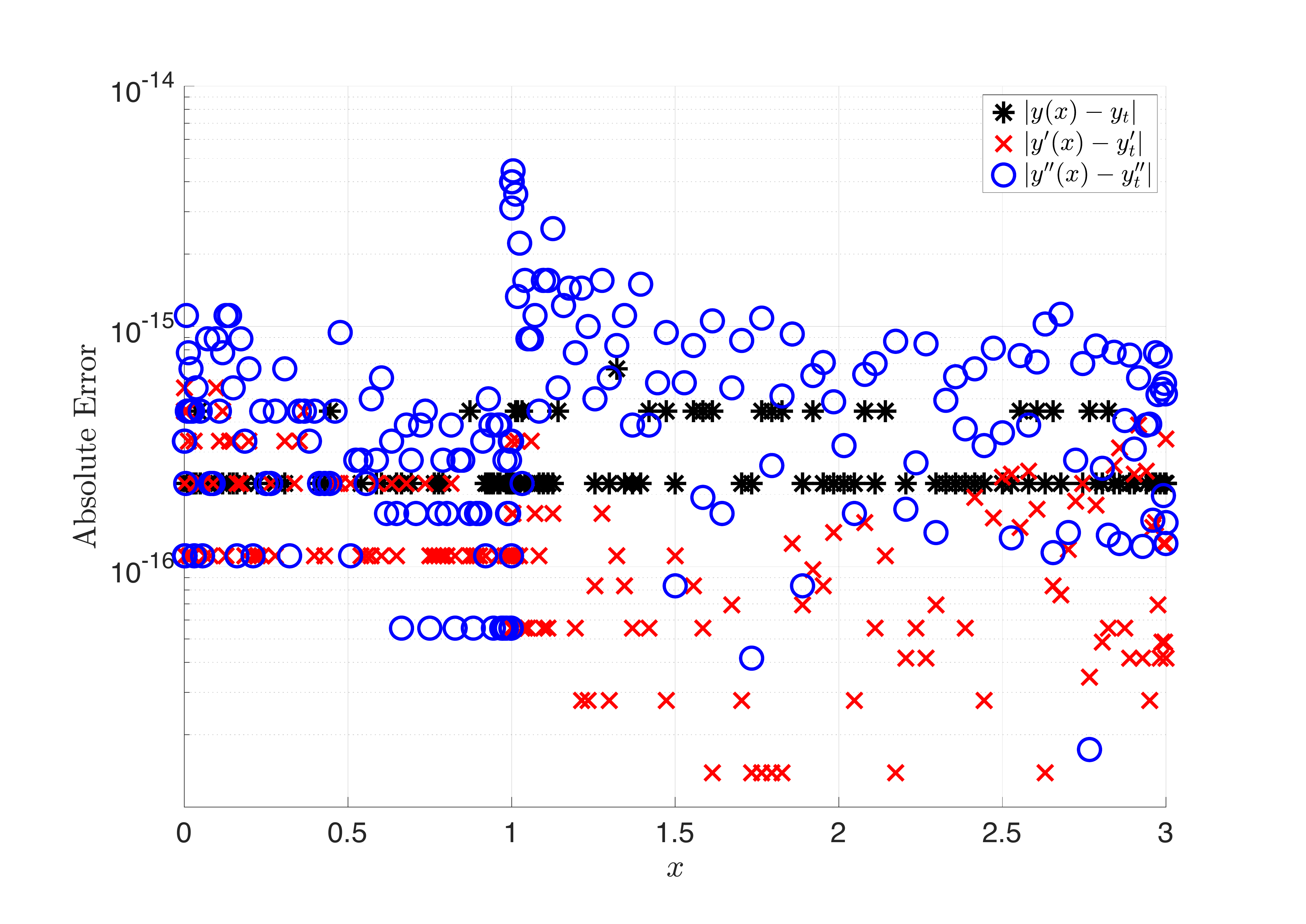}}%
\caption{Results of nonlinear-nonlinear differential equation using the TFC approach. It can be seen that the absolute error for the function and subsequent derivatives are around $10^{-15}$ to $10^{-16}$.}
\label{fig:ex3}
\end{figure}

The last example looks to complete the study with a sequence of two nonlinear differential equations. Consider the nonlinear-nonlinear DE sequence such that,
\begin{equation}
    y'' - ay'^2 = 0 \quad \text{subject to: } \begin{cases} y(0) = 2 \\ y(3) = 2 - \dfrac{\log(11,264)}{10} \end{cases} \text{and } \begin{cases} a = 1 \quad &\text{for } x \le 1 \\ a = 10 \quad &\text{for } x > 1 \end{cases}
\end{equation}
which admits the analytical solution of the form,
\begin{equation*}
    y(x) = \begin{cases} = 2 - \log(x+1) \quad &\text{for } x \le 1 \\ = 2 - \dfrac{9}{10}\log(2) - \dfrac{1}{10}\log(10x-8) \quad &\text{for } x > 1 \end{cases}
\end{equation*}
This example shares the same structure as that presented in Section \ref{sec:ex2} excepted the loss functions and Jacobians differ. For this differential equation the loss functions for the two segments are,
\begin{equation*}
    \mathcal{L}_{(1)}(x) = y''_{(1)}(x) - y'^2_{(1)}(x)
\end{equation*}
\begin{equation*}
    \mathcal{L}_{(2)}(x) = y''_{(2)}(x) - 10 \, y'^2_{(2)}(x)
\end{equation*}
The partials of the loss functions are of the form of Eq. (\ref{eq:Jacobian}) with the terms defined in \ref{sec:nonlinear_nonlinear} in Eqs. (\ref{eq:c1}-\ref{eq:c6}). The same initialization process as Section \ref{sec:ex2} was used to determine the initial guess for the unknown coefficients,
\begin{equation*}
    \Xi_0 = \begin{Bmatrix} \B{0}\T & 1.30685 & -0.5 & \B{0}\T \end{Bmatrix}\T.
\end{equation*}
Using $N = 100$ points and $m = 60$ basis functions per segment and the same convergence tolerance as the prior example with $\varepsilon = 10^{-15}$, the method reached machine error accuracy in $7$ iterations, only. The results of this test are displayed in Figs. \ref{fig:ex3:function} and \ref{fig:ex3:error} where it can be seen that the DE is solved with the same accuracy as the other cases presented. 

\section{Conclusions}

This paper explores the application of the \tfc\ \cite{U-ToC} to boundary-value problems on hybrid systems (physical systems defined by DE sequences). These types of problems were initially motivated by optimal control problems derived from the indirect method where the control undergoes discrete switches over the trajectory of the solution (i.e. fuel optimal landing, bang-bang control, etc.) and the system must satisfy boundary constraints (e.g., docking, landing).

TFC applied to hybrid systems looks to split the domain into multiple segments coinciding with the differential equations. Continuity is enforced between segments by relative parameters which are solved along with the unknown coefficients defining the free function. In general, the proposed method completely eliminates the need to use numerical shooting methods to ensure continuity since it is always satisfied through the use of TFC's \ces. Furthermore, the technique of this approach can be generalized for $n$ segments which admits a recursive form of the final matrix which is block diagonal.

In addition to the theoretical development, three numerical tests highlight the generality of this approach. The first numerical example produced a solution to a sequence composed of two linear differential equations. This was solved by linear least-squares requiring no iterations. The next two examples explore the linear-nonlinear and the nonlinear-nonlinear differential equations sequences and necessitate a nonlinear least-squares approach since at least one of the differential equations is nonlinear. The line connecting the two boundary values is adopted as the initial guess. This linear initial guess is easy to implement and it has always generated convergence in all the extensive hybrid systems tests performed to validate the proposed method.

The major result of all numerical tests is that the method solves the problems at machine level accuracy for all cases. The accuracy and speed results obtained are similar to those presented in past work on the \tfc\ \cite{LDE,NDE} applied to linear and nonlinear ODEs.

The main focus of this work was specific to $2^{nd}$-order differential equations where $C^1$ continuity is required. This research was initially motivated by the study of fuel-optimal planetary landing through the indirect method, which results in a $2^{nd}$-order system of BVPs. The use of this extension is necessary for the application of the Theory of Functional Connections to solve this problem and is an area of ongoing research \cite{fuel_opt}. Regardless, this method is general and can be used on any arbitrary order such that $C^{n-1}$ is satisfied for an $n^{th}$-order DE. 

\appendix
\section{Linear-Nonlinear Partials}\label{sec:linear_nonlinear}
The partial derivatives of the loss function of the linear-nonlinear differential sequence are provided below for completeness

\begin{align}
    \dfrac{\partial \mathcal{L}_{(1)}}{\partial \B{\xi}_{(1)}} &= \Big[\B{h}''_{(1)}(x) - \beta''_1(x)\B{h}_{(1)}(x_0) -\beta''_2(x) \B{h}_{(1)}(x_1) - \beta''_3(x)\B{h}'_{(1)}(x_1) \label{eq:b1}\\ & \qquad + \B{h}_{(1)}(x) - \beta_1(x)\B{h}_{(1)}(x_0) -\beta_2(x) \B{h}_{(1)}(x_1) - \beta_3(x)\B{h}'_{(1)}(x_1)\Big]\T \nonumber\\
    \dfrac{\partial \mathcal{L}_{(1)}}{\partial y_1} &= \beta''_2(x) + \beta_2(x) \label{eq:b2}\\
    \dfrac{\partial \mathcal{L}_{(1)}}{\partial y'_1} &= \beta''_3(x) + \beta_3(x) \label{eq:b3}\\
    \dfrac{\partial \mathcal{L}_{(2)}}{\partial \B{\xi}_{(2)}} &= \Big[\B{h}''_{(2)}(x) - \beta''_4(x)\B{h}_{(2)}(x_1) -\beta''_5(x) \B{h}'_{(2)}(x_1) - \beta''_6(x)\B{h}_{(2)}(x_f) \label{eq:b4}\\ & \qquad + y'_{(2)}(x) \Big(\B{h}_{(2)}(x) - \beta_4(x)\B{h}_{(2)}(x_1) -\beta_5(x) \B{h}'_{(2)}(x_1) - \beta_6(x)\B{h}_{(2)}(x_f) \Big) \nonumber\\ & \qquad + y_{(2)}(x) \Big(\B{h}'_{(2)}(x) - \beta'_4(x)\B{h}_{(2)}(x_1) - \beta'_5(x) \B{h}'_{(2)}(x_1) - \beta'_6(x)\B{h}_{(2)}(x_f) \Big)\Big]\T \nonumber\\
    \dfrac{\partial \mathcal{L}_{(2)}}{\partial y_1} &= \beta''_4(x) + y_{(2)}(x)\beta'_4(x) + y'_{(2)}(x)\beta_4(x) \label{eq:b5}\\
    \dfrac{\partial \mathcal{L}_{(2)}}{\partial y'_1} &= \beta''_5(x) + y_{(2)}(x)\beta'_5(x) + y'_{(2)}(x)\beta_5(x) \label{eq:b6}
\end{align}

\section{Nonlinear-Nonlinear Partials}\label{sec:nonlinear_nonlinear}
The partial derivatives of the loss function of the nonlinear-nonlinear differential sequence are provided below for completeness
\begin{align}
    \dfrac{\partial \mathcal{L}_{(1)}}{\partial \B{\xi}_{(1)}} &= \Big[\B{h}''_{(1)}(x) - \beta''_1(x)\B{h}_{(1)}(x_0) -\beta''_2(x) \B{h}_{(1)}(x_1) - \beta''_3(x)\B{h}'_{(1)}(x_1) \label{eq:c1} \\ & \qquad - 2 y'_{(1)}(x) \Big( \B{h}'_{(1)}(x) - \beta'_1(x)\B{h}_{(1)}(x_0) -\beta'_2(x) \B{h}_{(1)}(x_1) - \beta'_3(x)\B{h}'_{(1)}(x_1) \Big)\Big]\T \nonumber\\
    \dfrac{\partial \mathcal{L}_{(1)}}{\partial y_1} &= \beta''_2(x) - 2 y'_{(1)}(x) \beta'_2(x)\label{eq:c2}\\
    \dfrac{\partial \mathcal{L}_{(1)}}{\partial y'_1} &= \beta''_3(x) - 2 y'_{(1)}(x) \beta'_3(x) \label{eq:c3}\\
    \dfrac{\partial \mathcal{L}_{(2)}}{\partial \B{\xi}_{(2)}} &= \Big[\B{h}''_{(2)}(x) - \beta''_4(x)\B{h}_{(2)}(x_1) -\beta''_5(x) \B{h}'_{(2)}(x_1) - \beta''_6(x)\B{h}_{(2)}(x_f) \label{eq:c4}\\ & \qquad - 20 y'_{(2)}(x) \Big(\B{h}'_{(2)}(x) - \beta'_4(x)\B{h}_{(2)}(x_1) -\beta'_5(x) \B{h}'_{(2)}(x_1) - \beta'_6(x)\B{h}_{(2)}(x_f) \Big)\Big]\T \nonumber\\
    \dfrac{\partial \mathcal{L}_{(2)}}{\partial y_1} &= \beta''_4(x) - 20 y'_{(2)}(x) \beta'_4(x)\label{eq:c5}\\
    \dfrac{\partial \mathcal{L}_{(2)}}{\partial y'_1} &= \beta''_5(x) - 20 y'_{(2)}(x) \beta'_5(x)\label{eq:c6}
\end{align}

\section*{}
\noindent
\textbf{Conflict of Interest}: The authors declare no conflict of interest.

\section*{}
\noindent
\textbf{Funding}: This work was supported by a NASA Space Technology Research Fellowship, Johnston [NSTRF 2019] Grant \#: 80NSSC19K1149.

\bibliography{main}

\begin{thebibliography}{10}
\providecommand{\url}[1]{#1}
\csname url@samestyle\endcsname
\providecommand{\newblock}{\relax}
\providecommand{\bibinfo}[2]{#2}
\providecommand{\BIBentrySTDinterwordspacing}{\spaceskip=0pt\relax}
\providecommand{\BIBentryALTinterwordstretchfactor}{4}
\providecommand{\BIBentryALTinterwordspacing}{\spaceskip=\fontdimen2\font plus
\BIBentryALTinterwordstretchfactor\fontdimen3\font minus
  \fontdimen4\font\relax}
\providecommand{\BIBforeignlanguage}[2]{{%
\expandafter\ifx\csname l@#1\endcsname\relax
\typeout{** WARNING: IEEEtran.bst: No hyphenation pattern has been}%
\typeout{** loaded for the language `#1'. Using the pattern for}%
\typeout{** the default language instead.}%
\else
\language=\csname l@#1\endcsname
\fi
#2}}
\providecommand{\BIBdecl}{\relax}
\BIBdecl

\bibitem{VSS}
{Edwards, C., et al.(Eds.)}, \emph{Advances in Variable Structure and Sliding
  Mode Control}.\hskip 1em plus 0.5em minus 0.4em\relax Berlin, Heidelberg:
  Springer Berlin Heidelberg, 2006.

\bibitem{shooting_TPBVP}
\BIBentryALTinterwordspacing
D.~D. Morrison, J.~D. Riley, and J.~F. Zancanaro, ``Multiple shooting method
  for two-point boundary value problems,'' \emph{Commun. ACM}, vol.~5, no.~12,
  pp. 613--614, Dec. 1962. [Online]. Available:
  \url{http://doi.acm.org/10.1145/355580.369128}
\BIBentrySTDinterwordspacing

\bibitem{shooting_bang}
G.~J. Lastman, ``A shooting method for solving two-point boundary-value
  problems arising from non-singular bang-bang optimal control problems,''
  \emph{International Journal of Control}, vol.~27, no.~4, pp. 513--524, 1978.

\bibitem{shooting_BVP}
\BIBentryALTinterwordspacing
M.~Osborne, ``On shooting methods for boundary value problems,'' \emph{Journal
  of Mathematical Analysis and Applications}, vol.~27, no.~2, pp. 417 -- 433,
  1969. [Online]. Available:
  \url{http://www.sciencedirect.com/science/article/pii/0022247X69900596}
\BIBentrySTDinterwordspacing

\bibitem{LineShoot}
\BIBentryALTinterwordspacing
S.~M. Filipov, I.~D. Gospodinov, and I.~Faragó, ``Replacing the finite
  difference methods for nonlinear two-point boundary value problems by
  successive application of the linear shooting method,'' \emph{Journal of
  Computational and Applied Mathematics}, vol. 358, pp. 46 -- 60, 2019.
  [Online]. Available:
  \url{http://www.sciencedirect.com/science/article/pii/S0377042719301232}
\BIBentrySTDinterwordspacing

\bibitem{shooting_converg}
\BIBentryALTinterwordspacing
R.~Weiss, ``The convergence of shooting methods,'' \emph{BIT Numerical
  Mathematics}, vol.~13, no.~4, pp. 470--475, Dec 1973. [Online]. Available:
  \url{https://doi.org/10.1007/BF01933411}
\BIBentrySTDinterwordspacing

\bibitem{shooting_err1}
\BIBentryALTinterwordspacing
P.~Marzulli and G.~Gheri, ``Estimation of the global discretization error in
  shooting methods for linear boundary value problems,'' \emph{Journal of
  Computational and Applied Mathematics}, vol.~28, pp. 309 -- 314, 1989.
  [Online]. Available:
  \url{http://www.sciencedirect.com/science/article/pii/0377042789903427}
\BIBentrySTDinterwordspacing

\bibitem{shooting_err2}
\BIBentryALTinterwordspacing
P.~Marzulli, ``Global error estimates for the standard parallel shooting
  method,'' \emph{Journal of Computational and Applied Mathematics}, vol.~34,
  no.~2, pp. 233 -- 241, 1991. [Online]. Available:
  \url{http://www.sciencedirect.com/science/article/pii/037704279190045L}
\BIBentrySTDinterwordspacing

\bibitem{FEM}
\BIBentryALTinterwordspacing
J.~N. Reddy, ``{An Introduction to the Finite Element Method},'' \emph{Journal
  of Pressure Vessel Technology}, vol. 111, no.~3, pp. 348--349, 08 1989.
  [Online]. Available: \url{https://doi.org/10.1115/1.3265687}
\BIBentrySTDinterwordspacing

\bibitem{U-ToC}
\BIBentryALTinterwordspacing
D.~Mortari, ``{The Theory of Connections: Connecting Points},'' \emph{MDPI
  Mathematics}, vol. 5 (57), no.~57, pp. 1--15, 2017. [Online]. Available:
  \url{http://www.mdpi.com/2227-7390/5/4/57}
\BIBentrySTDinterwordspacing

\bibitem{LDE}
\BIBentryALTinterwordspacing
------, ``{Least-squares Solution of Linear Differential Equations},''
  \emph{MDPI Mathematics}, vol.~5, no.~48, pp. 1--18, 2017. [Online].
  Available: \url{http://www.mdpi.com/2227-7390/5/4/48}
\BIBentrySTDinterwordspacing

\bibitem{NDE}
D.~Mortari, H.~Johnston, and L.~Smith, ``{High Accurate Least-Squares Solutions
  of Nonlinear Differential Equations},'' \emph{Journal of Computational and
  Applied Mathematics}, vol. 352, pp. 293--307, January 8--12 2019.

\bibitem{M-ToC}
D.~Mortari and C.~Leake, ``{The Multivariate Theory of Connections},''
  \emph{MDPI Mathematics}, vol.~7, no.~3, p. 296, 2019.

\bibitem{constraints}
H.~Johnston and D.~Mortari, ``{Linear Differential Equations Subject to
  Relative, Integral, and Infinite Constraints},'' in \emph{2018 AAS/AIAA
  Astrodynamics Specialist Conference Snowbird, UT, August 19--23, 2018}.\hskip
  1em plus 0.5em minus 0.4em\relax AAS/AIAA, 2018.

\bibitem{TFC_selected}
\BIBentryALTinterwordspacing
H.~Johnston, C.~Leake, Y.~Efendiev, and D.~Mortari, ``Selected applications of
  the theory of connections: A technique for analytical constraint embedding,''
  \emph{Mathematics}, vol.~7, no.~6, 2019. [Online]. Available:
  \url{https://www.mdpi.com/2227-7390/7/6/537}
\BIBentrySTDinterwordspacing

\bibitem{Colloc}
C.~Lanczos, \emph{Applied Analysis}.\hskip 1em plus 0.5em minus 0.4em\relax New
  York: Dover Publications, Inc., 1957, p. 504.

\bibitem{ChebCol}
K.~Wright, ``{Chebyshev Collocation Methods for Ordinary Differential
  Equations},'' \emph{The Computer Journal}, vol.~6, no.~4, pp. 358--365,
  January 1964.

\bibitem{fuel_opt}
E.~Schiassi, R.~Furfaro, H.~Johnston, and D.~Mortari, ``{Fuel-efficient Powered
  Descent Guidance on Planetary Bodies via Theory of Functional Connection 1:
  Solution of the Equations of Motion},'' in \emph{2019 AAS/AIAA Astrodynamics
  Specialist Conference, Portland, ME, August 11-15, 2019}.\hskip 1em plus
  0.5em minus 0.4em\relax AAS/AIAA, 2019.

\end{thebibliography}
\bibliographystyle{IEEEtran}

\end{document}